\documentstyle[12pt,twoside,a4mlnew,math30,etamak,etahead,ready]{article}

\def\revmark{\relax}

\def\myleftmark{Br\"uning and Lesch: On the $\eta$--invariant}
\newcommand{\comment}[1]{\relax}

\def\sectone{2}
\def\secttwo{3}
\def\sectthree{4}
\def\endproof{\hfill\qed\par\vspace{1em}\bigbreak\noindent}
\newcommand{\mls}{\hspace*{-4pt}}
\begin{document}
\pagestyle{empty}
\renewcommand{\tilde}{\widetilde}
\def\VMark{\relax}

\def\G4-1.3{{\rm (2.3)}}

\renewcommand{\thefootnote}{\arabic{footnote}}
\begin{center}
\LARGE\bf On the eta--invariant of certain non--local
boundary value problems

\bigskip
\bigskip
\Large Jochen Br\"uning and Matthias Lesch

\bigskip\normalsize
Institut f\"ur Mathematik,
Humboldt--Universit\"at zu Berlin,
Unter den Linden 6,
D--10099 Berlin,
Germany,\\
e-mail: bruening@mathematik.hu-berlin.de, lesch@mathematik.hu-berlin.de\\
Internet: http://spectrum.mathematik.hu-berlin.de/$^\sim$lesch

\bigskip
\bigskip
Revised Version 14 Nov 1997
\end{center}

\bigskip\bigskip
\tableofcontents
\newpage\hphantom{++}\newpage
\pagestyle{headings}
\setcounter{page}{1}
\section{Introduction}\label{Sec0}

A very intriguing feature of elliptic operators on compact manifolds
without boundary \revmark is the locality of their indices. Specifically, if $M$ denotes a compact
Riemannian spin manifold, $S\to M$ a spinor bundle, $E\to M$ a hermitian
coefficient bundle with unitary connection, and $D^E$ the Dirac operator
on $M$ with coefficients in $E$ then, by the Atiyah--Singer theorem,
\begin{equation}
   \ind D_+^E=\int_M \hat A(M)\wedge {\rm ch} \, E.\mylabel{I1}
\end{equation}
Here $D_+^E$ arises from splitting $S\otimes E$ under the involution
induced by the complex volume element on $M$.

If $M$ decomposes along a compact hypersurface, $N$, as $M=M_1\cup M_2$,
with $\pl M_i=N$ for $i=1,2$, then one is lead to ask whether
the obvious decomposition of the right hand side in \myref{I1} 
corresponds to a decomposition of the (essentially) self--adjoint operator
$D^E$ into self--adjoint operators $D_i^E$, defined in $M_i$ by suitable
boundary conditions on $N$, such that
\begin{equation}
  \ind D_{1,+}^E+\ind D_{2,+}^E=\ind D_+^E.\mylabel{I2}
\end{equation}
This question was answered in the affirmative by Atiyah, Patodi, and Singer
\cite{APS} who formulated the correct boundary conditions
(cf. Sec. \mls\sectone\ for details). More importantly, the resulting index formula
\myref{G4-1.6} displayed a new spectral invariant of self--adjoint
elliptic operators (defined on $N$) which they called the
$\eta$--invariant. It is not locally computable by a formula
as in \myref{I1} as can be seen from its behaviour under coverings. Nevertheless,
one can ask how the $\eta$--invariant behaves under splitting $N$
as $N_1\cup N_2$, and this is the problem we address in this work.

One motivation for posing this question may be seen in trying
to understand the signature theorem on manifolds with corners.
From a systematical point of view, splitting formulas for spectral
invariants should also be very useful for
computational purposes -- as illustrated nicely by the analytic torsion,
cf. \cite{Cheeger,Mullertorsion} -- and as a possible source of new invariants. Another recent motivation
is provided by topological quantum field theory.

The ''gluing law'' for $\eta$--invariants we prove here 
(Thms. \plref{Srev-glu5}, \plref{Srev-glu6})
is not new; cf. Sec.\sectone\ for an account of previous work. Our proof, however,
attacks the problem directly  on the cut manifold, $M^{\rm cut}$, by 
analizing families of ''generalized Atiyah--Patodi--Singer boundary
value problems.'' 
These new abstract boundary conditions are defined
by three simple axioms
(\myref{G5-1.31a}--\myref{G5-2.25} below) which are designed
in such a way that the heat kernel of the model operator
is explicitly computable. Incidentally, our formula generalizes
a result of Sommerfeld in the scalar case. Moreover, under this
class we find the spectral boundary conditions introduced by
Atiyah, Patodi, and Singer as well as the (local) absolute
and relative boundary conditions for the Gau{\ss}--Bonnet operator. Thus, our
method gives a uniform way to derive the asymptotic expansion of the
heat trace in both cases, generalizing in particular recent work by
Grubb and Seeley \cite{GrubbSeeley} (cf. Thm. \mls\plref{S1-2.7}).
The family we define
interpolates between the ''uncut manifold'' (the case of smooth
transmission) and actual Atiyah--Patodi--Singer boundary value problems;
this is similar to Vishik's approach to the splitting behavior of
the analytic torsion, and we hope to exploit this further in another
publication. The special structure of our family, on the other hand, 
resembles closely the finite--dimensional variations constructed by Lesch
and Wojciechowski \cite{LeschW}. This allows us to produce
explicit variation formulas (Thm. \mls\plref{S6-2.5}). We evaluate them
using the vanishing of the noncommutative residue on pseudodifferential
idempotents and a special symmetry of the cutting problem.

\pagebreak[3]
The plan of the paper is as follows: In Section \sectone, we review
some abstract facts on $\eta$--invariants and previous work
on the gluing law. All results are presented in Section \secttwo\
while the details of most proofs are carried out in Section
\sectthree.

This work was supported by Deutsche Forschungsgemeinschaft and the
GADGET network of the EU.

\section{Generalities}\label{Sec1}

In this section we briefly review some more or less well known properties
of $\eta$--invariants which are needed below, together with some of the
previous work leading to the gluing law.

The $\eta$--invariant was introduced in the seminal work
\cite{APS} by Atiyah, Patodi, and Singer. They considered
the signature operator on a smooth oriented Riemannian
manifold, $M$, with compact boundary $\pl M=N, \dim M=m=4k$.\revmark\ 
The signature operator is the operator $D=d+\delta$ restricted to
the space of self--dual forms (cf. \myref{rev-G2.6} below).
Assuming that the metric is a product in a neighborhood 
\alpheqn
\begin{equation}
U\simeq[0,1)\times N\mylabel{G4-1.1a}
\end{equation} of the boundary, separation of
variables leads to the representation
\begin{equation}
D=\gamma(\frac{\pl}{\pl x}+A).
  \mylabel{G1-1.2}
\end{equation}
Here, we use the decomposition of a smooth form, $\alpha$, as
$\alpha=dx\wedge \ga_1(x)+\ga_2(x)$. Thus, the operator on the right acts on $\cinfz{(0,1),
\Omega(N)\oplus\Omega(N)}$, $\Omega(N)$ the smooth forms on $N$, and one 
has\revmark
\begin{equation}
\gamma=\mat{0}{1}{-1}{0}\otimes I,\quad A=\mat{0}{-1}{-1}{0}\otimes
   (d_N+\delta_N).\mylabel{Grev-2.1c}
\end{equation}   
\reseteqn   
Thus $A$ is symmetric, and we have the relations
\begin{equation}
\gamma^2=-I,\quad \gamma^*=-\gamma, \quad \gamma A+A\gamma=0.
  \mylabel{G1-1.3}
\end{equation}
A symmetric operator of type \myref{G1-1.2}
does not in general admit local
boundary conditions
which define a self--adjoint extension (cf., however,
\cite{GilkeySmith} and \cite{Singer}), even though local boundary
conditions do exist in the special case \myref{Grev-2.1c} i.e.
the absolute and relative boundary conditions. 
But there is always a nonlocal
boundary condition given (essentially) by the Calder\'on projector
\cite{Calderon}. Thus we introduce the boundary condition
\alpheqn
\begin{equation}
   P_{>0}(A) u(0)=0,
   \mylabel{G4-1.3a}
\end{equation}
where $P_{>0}(A)$ is the orthogonal projection onto the subspace
spanned by eigenvectors of $A$ with positive eigenvalues. To define
a symmetric operator, this
needs to be supplemented by
\begin{equation}
  P_\sigma u(0)=0,
  \mylabel{G4-1.3b}
\end{equation}
\reseteqn
where $P_\sigma$ projects onto a Lagrangian subspace of $\ker A$
with respect to the symplectic form (note that $\dim\ker A$
is even)
$$\omega(u,v):=<\gamma u,v>,\quad u,v\in\ker A,$$
and such a space can always be viewed as the $+1$--eigenspace
of an involution, $\sigma$, on $\ker A$ satisfying
\alpheqn
\begin{equation}
\sigma\gamma+\gamma\sigma=0;\mylabel{G4-1.4a}
\end{equation}
then
\begin{equation}
  P_\sigma=\frac 12 (I+\sigma).\mylabel{G4-1.4b}
\end{equation}
In the case at hand, a convenient choice of $\sigma$ is
(Clifford multiplication by) the
complex volume element, $\omega_M$, i.e. we put
$\sigma_0:=\omega_M|\ker A$ and observe that it takes the form
$$\omega_M=\mat{0}{1}{1}{0}\otimes \omega_N,$$
where $\omega_N$ denotes the complex volume element on $N$.
\reseteqn

It is
not hard to see that these data define a self--adjoint extension
of $D$, $D_{\sigma_0}$, which anticommutes with $\omega_M$. Then
the signature operator, $D_S$, for a manifold with boundary
is the closure of
$$D_{\sigma_0}\big|\cd(D_{\sigma_0})\cap\left\{u\in \Omega(\ovl{M})\;|\;
    \omega_M u=u\right\},$$
and \cite[Thm. (I.3.10)]{APS} asserts that
$D_S$ is a Fredholm operator with
\begin{equation}
   \ind D_S=\int_M L(M)-\frac 12(\eta(B)+\dim\ker B).
   \mylabel{G4-1.5}
\end{equation}
Here, $L(M)$ denotes the Hirzebruch $L$--form and the operator $B$
is defined by a representation of $D_S$ in $U$ analogous to
\myref{G1-1.2}. In fact, near $\pl M$ we have
\begin{eqnarray*}
    D_S&=&  \omega_N\left(\pl_x+\omega_N(d_N+\delta_N)\right)\\
       &=:&  \omega_N\left(\pl_x+ B\right),
\end{eqnarray*}
and a core is given by the space (with obvious notation)\revmark
\begin{equation}
\cd(D_S)=\left\{u\in \Omega(\ovl{M})\;|\; P_{\ge 0}(B) u(0)=0,
   \;\omega_M u=u\right\}.
  \mylabel{rev-G2.6}
\end{equation}

Rewriting \myref{G4-1.5} in terms of the signature of $M$
(as a manifold with boundary) gives \cite[Thm. (I.4.14)]{APS}
\begin{equation}
   \sign M=\int_M L(M)-\frac 12 \eta(B),
   \mylabel{G4-1.6}
\end{equation}
and thus an analytic interpretation of the additivity of the signature
under cutting along a separating hypersurface.

The $\eta$--invariant figuring in \myref{G4-1.5} and \myref{G4-1.6}
is derived from a meromorphic function generalizing the $\zeta$--function
of an elliptic operator. It is convenient to derive the main properties
of these functions in an abstract functional analytic setting. Thus
consider a self--adjoint
operator, $A$, with dense domain, $\cd(A)$, in some Hilbert space,
$H$. If we assume that
\begin{equation}
   (A+i)^{-1}\in C_p(H), \quad\mbox{\rm for some}\; p>0,
   \mylabel{G4-1.7}
\end{equation}
(where $C_p$ denotes the Schatten--von Neumann class of order $p$) then
the function
\begin{equation}
     \eta(A;s):= \frac{1}{\Gamma(\frac{s+1}{2})}\int_0^\infty
     t^{(s-1)/2} \tr_H\big(A e^{-t A^2}\big)dt=
     \sum_{\gl\in\spec A\setminus\{0\}} (\sgn \gl) |\gl|^{-s}
     \mylabel{G4-1.8}
\end{equation}
is holomorphic for large $\Re s$. More generally, if
$B:\cd(A)\to H$ is any bounded operator satisfying
\begin{equation}
    P_0(A)BP_0(A)=0, \mylabel{G4-1.9}
\end{equation}
$P_0(A)$ the orthogonal projection onto $\ker A$, then the same is true of
\begin{eqnarray}
   \eta(A,B;s)&:=& \frac{1}{\Gamma(\frac{s+1}{2})}\int_0^\infty
     t^{(s-1)/2} \tr_H\big(B e^{-t A^2}\big)dt\nonumber\\
     &=&\sum_{\gl\in\spec A\setminus\{0\}} (\tr_{\ker(A-\gl)} B) |\gl|^{-s-1}.
     \mylabel{G4-1.10}
\end{eqnarray}
Here, by slight abuse of notation,
$\tr_{\ker(A-\gl)} B:=\tr(P_\gl B)$ where
$P_\gl$ is the orthogonal projection onto the $\gl$--eigenspace of $A$.
\revmark\ 
It is very important to determine conditions on $A$ and $B$ which
guarantee the existence of a meromorphic extension of \myref{G4-1.10} to
the whole complex plane. The standard source of such an extension
is an asymptotic expansion
\begin{equation}
    \tr_H\big(B e^{-t A^2}\big)\sim_{t\to 0+}
    \sum_{\begin{array}{c} \SST\Re \ga\to \infty\\
         \SST 0\le k\le k(\ga)\end{array}}
    a_{\ga k}(A,B)\, t^\ga\, \log^kt.
   \mylabel{G4-1.11}
\end{equation}
The notation used means, of course, that
$\{\ga\in\C\,|\, a_{\ga k}(A,B)\not=0\;\mbox{\rm for some}\; k\in \Z_+, k\le k(\ga)\}$
is a countable subset of $\C$ whose real parts accumulate at most at $\infty$.

Using the notation $f(s)=:\DST\sum_k\TST \Res_k f(s_0) (s-s_0)^{-k}$,
introduced in \cite{BS1}
for Laurent expansions, one has

\begin{lemma}{S4-1.1}
Under the conditions \myref{G4-1.7}, \myref{G4-1.9}, and \myref{G4-1.11},
$\eta$ extends to a meromorphic function on $\C$.

The poles are situated at the points $s_\ga=-2\ga-1$ and the principal
part of $\eta$ at $s_\ga$ is given by
$$\frac{1}{\Gamma(\frac{s+1}{2})}\sum_{k=0}^{k(\ga)}   a_{\ga,k}(A,B)
 (-1)^k k! 2^{k+1}(s-s_\ga)^{-k-1}.$$
In particular, the poles are of order
\renewcommand{\labelenumi}{{\rm (\arabic{enumi})}}
\begin{enumerate}
\item $k(\ga)+1$, if $\ga\not\in\Z_+$, and
\alpheqn
\begin{equation}
\Res_{k(\ga)+1} \eta(A,B;s_\ga)=\frac{(-1)^{k(\ga)}k(\ga)! 2^{k(\ga)+1}}%
   {\Gamma(-\ga)} a_{\ga,k(\ga)}(A,B),
  \mylabel{G4-1.11a}
\end{equation}
\item[and]
\item $k(\ga)$, if $\ga\in\Z_+$, and
\begin{equation}
\Res_{k(\ga)} \eta(A,B;s_\ga)=(-1)^{k(\ga)+\ga} \ga! k(\ga)! 2^{k(\ga)}
   a_{\ga,k(\ga)}(A,B).
   \mylabel{G4-1.11b}
\end{equation}
\reseteqn
\end{enumerate}
\end{lemma}

\pagebreak[3]

\begin{lemma}{S6-1.2.5} Under the conditions
\myref{G4-1.7} and \myref{G4-1.9} the following statements
are equivalent:
\begin{itemize}
\item[{\rm (i)}] $\tr_H\big(Be^{-tA^2}\big)$ has an asymptotic
expansion of type \myref{G4-1.11} which can be differentiated, i.e.
for $N,K>0$ we have
\begin{equation}\mylabel{G7-1.12}
   \Big|\pl_t^N\Big(\tr_H\big(Be^{-tA^2}\big)
      -\sum_{\begin{array}{c} \SST\Re \ga\le N+K\\
         \SST 0\le k\le k(\ga)\end{array}}
    a_{\ga k}(A,B)\, t^\ga\, \log^kt\Big)\Big|
      \le C_{N,K} t^K,\quad t\to 0.
\end{equation}
\item[{\rm (ii)}] $\Gamma(\frac{s+1}{2})\eta(A,B;s)$ is holomorphic
  in the half plane $\{s\in\C\,|\, \Re s> p\}$ and extends meromorphically
  to $\C$. Moreover, for $a,b\in\R$
  there exists  $s_0=s_0(a,b)>0$ such that $\Gamma(\frac{s+1}{2})\eta(A,B;s)$
  is holomorphic for $a\le \Re s\le b,|s|\ge s_0$ with estimate
\begin{equation}\mylabel{G7-1.13}
    \left|\Gamma({\TST\frac{s+1}{2}})\eta(A,B;s)\right|
    \le C(a,b,N) |s|^{-N},\quad a\le \Re s\le b,|s|\ge s_0,
\end{equation}
for any $N>0$.
\end{itemize}
\end{lemma}
\proof {(i)$\Rightarrow$(ii):}
In view of \myref{G4-1.7} and \myref{G4-1.9} $\Gamma(\frac{s+1}{2})\eta(A,B;s)$
is holomorphic in the half plane $\{ s\in\C\,|\, \Re s>p\}$ and extends
meromorphically to $\C$, by Lemma \plref{S4-1.1}. Integration by parts gives
\begin{equation}\mylabel{G7-1.14}
   \Gamma({\TST\frac{s+1}{2}})\eta(A,B;s)={\TST
    \frac{(-1)^N 2^N}{(s+1)(s+3)\cdot\ldots\cdot(s+2N-1)}}
    \int_0^\infty t^{(s-1)/2+N} \pl_t^N
    \tr_H\big(Be^{-tA^2}\big)dt.
\end{equation}
In view of \myref{G4-1.9} we have for $a\le \Re s\le b$
\begin{equation}\mylabel{G7-1.15}
  \left|\int_1^\infty t^{(s-1)/2+N} \pl_t^N
    \tr_H\big(Be^{-tA^2}\big)dt\right|\le
     C\int_1^\infty t^{(b-1)/2+N} e^{-\eps t}dt=:C_{N,b}.
\end{equation}
Furthermore, choosing $K$ such that $(a-1)/2+K+N>-1$, we may write
\begin{equation}\begin{array}{ll}
&\DST\int_0^1 t^{(s-1)/2+N} \pl_t^N\tr_H\big(Be^{-tA^2}\big)dt\\[1em]
=:&\DST\int_0^1 t^{(s-1)/2+N} \varphi_{K,N}(t) dt
       +\hspace*{-0.5em}\sum_{\begin{array}{c} \SST\Re \ga\le N+K\\
         \SST 0\le k\le k(\ga)\end{array}}\hspace*{-0.5em}
    a_{\ga k}(A,B)\int_0^1 t^{(s-1)/2+N} \pl_t^N t^\ga\, \log^kt dt
    \end{array}
    \mylabel{G7-1.16}
\end{equation}
with $|\varphi_{K,N}(t)|\le C_{K,N} t^K$. Hence, we have for
$a\le \Re s\le b$
\begin{equation}\mylabel{G7-1.17}
   \left| \int_0^1 t^{(s-1)/2+N} \varphi_{K,N}(t) dt\right|\le C_{N,K}.
\end{equation}
Using $\pl_t^N t^\ga\, \log^kt=
\sum\limits_{i=0}^k c_i t^{\ga-N}\log^i t$ we get

\begin{equation}\mylabel{G7-1.18}
   \int_0^1 t^{(s-1)/2+N}\pl_t^N t^\ga\, \log^kt=
   \sum_{i=0}^k c_i (-1)^i i! ((s+1)/2+\ga)^{-i-1}.
\end{equation}
Combining \myref{G7-1.14} through \myref{G7-1.18} we reach the conclusion.

\medskip
{(ii)$\Rightarrow$(i):} In view of the estimate \myref{G7-1.13}
we can apply the inverse Mellin transform to find, for $c>p$, 
$$\tr_H\big(Be^{-tA^2}\big)=\frac{1}{4\pi i}\int_{\Re s=c}
    t^{-(s+1)/2}\, \Gamma({\TST\frac{s+1}{2}})\eta(A,B;s) ds.$$
Moreover, we can shift the contour of integration to the left
and apply the Residue Theorem to get
$$\tr_H\big(Be^{-tA^2}\big)\sim_{t\to 0+}
   \frac 12\sum_{s\in\C}
   \Res_1\Big(t^{-(s+1)/2}\, \Gamma({\TST\frac{s+1}{2}})\eta(A,B;s)\Big).$$
Clearly, this asymptotic expansion can be differentiated.\endproof

\vspace{-1em}
\begin{absatz}{Remarks} 1)
Of course, $B:=I-P_0(A)$ gives the $\zeta$--function of $A^2$, 
$$\zeta_{A^2}({\TST\frac{s+1}{2}})=\eta(A,I-P_0(A);s).$$
In particular, we can read off the regularity at $0$ of
$\zeta_{A^2}$ provided that the asymptotic expansion of
$\tr_H\big(e^{-tA^2}\big)$ exists and does not contain
contributions to $\log^k t, k\in\N$.

2) If $A$ and $B$ are classical pseudodifferential
operators on a compact manifold, $M$, $\dim M=:m$, and
$A$ is self--adjoint elliptic of positive order, then \myref{G4-1.7}\revmark\
holds and we have an asymptotic expansion
\cite[Theorem 2.7]{GrubbSeeley}
\begin{equation}
  \tr_H\big(B e^{-tA^2}\big)\sim_{t\to 0+}
    \sum_{j=0}^\infty a_j(A,B)\, t^{(j-m-b)/2a}
      +\sum_{j=0}^\infty (b_j(A,B)\, \log t+c_j(A,B)) t^j
      \mylabel{G7-1.19}
\end{equation}
where $a:=\ord A, b:=\ord B$. Moreover, this asymptotic expansion
can be differentiated in view of the identity
$$\pl_t^N\tr_H\big(B e^{-tA^2}\big)=(-1)^N\tr_H\big(
   BA^{2N} e^{-tA^2}\big).$$
If, in addition, \myref{G4-1.9} holds then we can apply Lemma
\plref{S6-1.2.5} to conclude that \myref{G7-1.13} holds for
$A$ and $B$.

Note that in view of \myref{G7-1.19} and Lemma \plref{S4-1.1}, in this case
$\eta(A,B;s)$ has a meromorphic continuation to $\C$ with
{\em simple} poles.

The estimate \myref{G7-1.13} suffices to shift the contour of
integration and to deduce a short time asymptotic expansion. However,
for some classical pseudodifferential operators $A,B$ an even
stronger result holds: Namely, if $A$ has scalar principal
symbol then it follows from \cite{DG} that
$\eta(A,B;s)$ is of polynomial growth on finite vertical strips.
Since $\Gamma(\frac{s+1}{2})$ decays exponentially on finite
vertical strips this implies the estimate
\myref{G7-1.13}. However, our method of proving \myref{G7-1.13}
is completely elementary while \cite{DG} uses the machinery
of Fourier integral operators.
\end{absatz}

Given these preparations we define, under the assumptions of Lemma
\plref{S4-1.1} (actually, a partial expansion in \myref{G4-1.11}
would suffice), the {\it $\eta$--invariant of $A$} as
\alpheqn
\begin{equation}
  \eta(A):= \Res_0 \eta(A;0),
  \mylabel{G4-1.12}
\end{equation}
and, in view of the index formula \myref{G4-1.5}, the {\it
reduced $\eta$--invariant of } $A$ as
\begin{equation}
  \xi(A):=\frac 12 \left(\eta(A)+\dim \ker A\right).
  \mylabel{G4-1.12b}
\end{equation}
\reseteqn
Generally, $\eta(A)$ is difficult to compute. It is thus of great
importance that suitable one--parameter variations turn out to be
"locally computable" in the sense of
asymptotic expansions of the type \myref{G4-1.11}.

To deal with variations in the abstract framework above we now impose the following
assumptions. Consider a connected open subset, $J$, of $\R$ and for $a\in J$ a
family
\alpheqn
\begin{equation}
  A(a):\cd\longrightarrow H, \mylabel{G4-1.13a}
\end{equation}
of self--adjoint operators with fixed domain $\cd$, satisfying
\myref{G4-1.7}. 

Moreover, assume that this family has kernel of constant rank, i.e.
for $P_0(a):=P_0(A(a))$ we have
\begin{equation}
 \dim P_0(a)\quad\mbox{\rm is constant in $J$}.\mylabel{G4-1.13b}
\end{equation}

Likewise, let
\begin{equation}
  B(a):\cd\longrightarrow H, \mylabel{G4-1.13c}
\end{equation}
be another family of bounded operators satisfying \myref{G4-1.9}
which, in addition, commutes with $A(a)^2$ in the sense that
\begin{equation}
 \big[B(a),(A(a)^2-\zeta)^{-1}\big]=0,\quad a\in J,\quad \zeta\not\in \spec A(a)^2.
  \mylabel{G4-1.13d}
\end{equation}
Note that \myref{G4-1.9} and \myref{G4-1.13d} imply that\revmark
$$ B(a)=(I-P_0(a))B(a)(I-P_0(a)).$$
Finally, we assume that
\setlength{\blockwidth}{\textwidth}
\setlength{\blockmargin}{8mm}
\setlength{\blockrightmargin}{1.7cm}
\addtolength{\blockwidth}{-\blockmargin}
\addtolength{\blockwidth}{-\blockrightmargin}
\begin{numabsatz}
  \myitem{the families $(A(a))_{a\in J}, (B(a))_{a\in J}\subset\cl(\cd,H)$
   are strongly differentiable in $J$, with strongly continuous derivative.}
\end{numabsatz}
\reseteqn
Under these assumptions, the operator families $P_0(a)$ and
\begin{equation}
   \tilde A(a):=(I-P_0(a))A(a)+P_0(a)
   \mylabel{G4-1.14}
\end{equation}
are strongly differentiable, too. Using the representation
$$e^{-tA(a)^2}=\frac{(m-1)!\; t^{1-m}}{2\pi i}\int_\Gamma e^{-t\zeta}
    (A(a)^2-\zeta)^{-m} d\zeta,$$
with $\Gamma$ a suitable contour, one can easily derive the identity
\begin{eqnarray*}
    \frac{\pl}{\pl a} \tr_H \left[ B(a) e^{-t A(a)^2}\right]
    &=& \tr_H\left[ B'(a) e^{-tA(a)^2}\right] \\
    &&+ t \frac{\pl}{\pl t} \tr_H\left[ B(a) \left(\frac{d}{da} A(a)^2\right)
      \tilde A(a)^{-2} e^{-t A(a)^2}\right].
\end{eqnarray*}
Our assumptions imply the absolute and locally uniform convergence
of the relevant $t$--integrals, and we arrive at
\begin{lemma}{S4-1.2} Under the assumptions \myref{G4-1.7} and 
{\rm (\plref{G4-1.13a}-e)} we have
the identity
\begin{eqnarray}\frac{\pl}{\pl a} \eta(A(a),B(a);s)&=&
    \eta(A(a),B'(a);s)\nonumber\\
  &&-\frac{s+1}{2} \eta(A(a), 
     B(a) \left(\frac{d}{da} A(a)^2\right)
      \tilde A(a)^{-2};s).
   \mylabel{G4-1.15}
\end{eqnarray}
\end{lemma}

\vspace{1em}\noindent
If we assume in addition that
\def\schnoerkel{{\rm d'}}
\addtocounter{equation}{-3}
\begin{equation}
 [B(a),(A(a)-\zeta)^{-1}]=0\quad \mbox{for}\quad a\in J, \zeta\not\in\spec
   A(a),
\end{equation}
\resetschnoerkel
\addtocounter{equation}{2}
then \myref{G4-1.15} simplifies to
\begin{equation}\frac{\pl}{\pl a} \eta(A(a),B(a);s)=
    \eta(A(a),B'(a);s)-(s+1)\; \eta(A(a),A'BA\tilde A(a)^{-2};s).
   \mylabel{G4-1.16}
\end{equation}
So, if both sides extend meromorphically to $\C$ then
\myref{G4-1.16} holds in $\C$, too.
We note in particular that
\begin{equation}\frac{\pl}{\pl a} \eta(A(a);s)=
    -s\;\eta(A(a),A'(a);s).
   \mylabel{G4-1.17}
\end{equation}
Thus we obtain the well known
\begin{cor}{S4-1.3} Assume \myref{G4-1.7}, {\rm (\plref{G4-1.13a},b,e),} and \myref{G4-1.11}
with $A(a)$ and $A'(a)$ in place of $B$. Then, for $k\in \Z_+$,
\begin{eqnarray}
    \frac{d}{da} \Res_{k} \eta(A(a);0)&=&
   -\Res_{k+1} \eta(A(a),A'(a);0)\nonumber\\
   & =&\frac{(-1)^{k+1} k! 2^{k+1}}{\sqrt{\pi}}a_{-1/2,k} (A(a),A'(a)).
   \mylabel{G4-1.18}
\end{eqnarray}
\end{cor}

\vspace{1em}\noindent
The condition \myref{G4-1.13b} is not satisfied in interesting situations.
One can get rid of it in choosing a real number
$c>0$ so that $c\not\in \spec(A(a))$ for $a$ near $a_0\in J$. Then we put
$\tilde P_{<c}(a):=P_{<c}(a)P_{>-c}(a),
\tilde P_{>c}(a):=I-\tilde P_{<c}(a)$ and
replace $A(a)$ by $A^c(a):=\tilde P_{>c}(a) A(a) + \tilde P_{<c}(a)$
and $B(a)$ by $B^c(a):=\tilde P_{>c}(a)B(a)\tilde P_{>c}(a) + \tilde
P_{<c}(a)$, obtaining the
modified $\eta$--function $\eta^c(A(a),B(a);s):=\eta(A^c(a),B^c(a);s)$.
$\eta^c$ admits, near $a_0$, the same analysis as outlined for $\eta$
with \myref{G4-1.13b}, and from \myref{G4-1.10} we obtain
\begin{equation}
   (\eta-\eta^c)(A(a),B(a);s)=\sum_{\begin{array}{c}
      \SST\gl\in\spec A(a)\\
      \SST 0<|\gl|< c \end{array}}
      |\gl|^{-s-1} \tr_{\ker (A(a)-\gl)} B(a)-\dim \tilde P_{<c}(a).
      \mylabel{G4-1.19}
\end{equation}

This is a smooth function of $a$ and holomorphic in $s\in\C$;
on the other hand, the negative $t$-powers in the expansion \myref{G4-1.11}
are unaffected if we modify $A$ and $B$ by an operator of finite
rank.
Evaluating \myref{G4-1.19} with $B(a):=A(a)$ we obtain
$$\frac 12 (\eta-\eta^c)(A(a);s)+\frac 12 \dim \ker A(a)=
    \sum_{\begin{array}{c}
      \SST\gl\in\spec A(a)\\
      \SST 0<|\gl|< c \end{array}}
      \frac 12 (\sgn \gl |\gl|^{-s}-1),$$
and consequently

\begin{lemma}{S4-1.3.5} Assume 
that the family $A(a)_{a\in J}$ satisfies 
\myref{G4-1.7}, {\rm (\plref{G4-1.13a},e),} and \myref{G4-1.11}
with $A(a)$ and $A'(a)$ in place of $B$. Then, for $a, a_0\in J$,
\begin{equation}
       \xi(A(a))-\xi(A(a_0))+\frac{1}{\sqrt{\pi}}\int_{a_0}^a
       a_{-1/2,0}(A(a),A'(a))da\;\in\Z.
   \mylabel{G4-1.20}
\end{equation}
\end{lemma}

This implies that the function
\begin{equation}
   \tau(A(a)):= e^{2\pi i \xi(A(a))}\mylabel{G4-1.21}
\end{equation}
is always smooth in $a\in J$ under our assumptions; the invariant
$\tau$ was introduced in \cite{DF}.

Instead of $\tau$ we will use the {\it reduced $\eta$--invariant}
\revmark
\begin{equation}
    \ovl{\eta}(A(a)):=\xi(A(a))\,\mod\,\Z,
\end{equation}
i.e. the image of the 
$\tau$--invariant under the diffeomorphism $S^1\cong \R/\Z$.

If the asymptotic expansion of $\tr_H\big[A^{j}(a) e^{-t A(a)^2}\big]$
does not contain terms of the form $t^\ga\, \log^k t$ with $\ga<0$
and $k\in\N$ for $j=0,1$ -- as it is the case for (classical) elliptic
pseudodifferential operators on compact manifolds, cf.
the remarks after Lemma \plref{S6-1.2.5} --
then it follows from Lemmas \plref{S4-1.1} and \plref{S4-1.2} that
$0$ is at most a simple pole of $\eta$ and that the residue is a homotopy
invariant. This is the basis for proving that $\eta(A;s)$ is, in fact,
regular at $s=0$ if $A$ happens to be a (classical) pseudodifferential
operator on a compact manifold, cf. \cite[Sec. 3.8]{Gi}.
More generally, Wodzicki \cite{Wod1,Wod2}\revmark\  observed the remarkable fact that, in this
class of operators,
\begin{equation}
\res B:= ({\rm ord}\,A) \Res_1 \eta(A,B;-1)=-2 ({\rm ord}\, A)\, a_{0,1}(A,B)
   \mylabel{G5-1.23}
\end{equation}   
defines the unique trace (up to a constant) on classical
pseudodifferential operators
if $A$ is elliptic of positive order, ${\rm ord}\,A$. 
Wodzicki also observed the following result, which
is stated without proof in his thesis \cite{Wod1}:\revmark

\begin{lemma}{S4-1.4} If $B$ is a classical pseudodifferential
operator on a compact manifold and an idempotent, then
$$\res B=0.$$
\end{lemma}

\noindent
The only proof we know of shows that the statement of this
lemma follows from the regularity at $0$
of the $\eta$--function for
general classical elliptic pseudodifferential operators on a
compact manifold.
For completeness we indicate that these facts are actually
equivalent.

\begin{lemma}{S6-1.5} The assertion of Lemma \plref{S4-1.4} is equivalent
to the following:
Let $P$ be a self--adjoint classical elliptic pseudodifferential operator
of positive order on the compact manifold $M$. Then 
$$\Res_1\eta(P;0)=0.$$
\end{lemma}
\proof
1. First we assume Lemma \plref{S4-1.4}. Let $P$ be a self--adjoint classical
elliptic pseudodifferential operator of order $d$ on a compact manifold, $M$.
We consider the pseudo\-dif\-ferential operator
$$\sgn P:=P |P|^{-1}:x\mapsto \casetwo{|P|^{-1}Px}{x\in\ker P^\perp,}{0}{x\in\ker P.}$$
We find
$$\eta(P^2,\sgn P;s)=\sum_{\gl\in\spec P} (\sgn \gl) |\gl|^{-s-1}=\eta(P;s+1)$$
an hence in view of \myref{G5-1.23} 
\begin{equation}
0=\res \sgn P= ({\rm ord}\,P)\Res_1\eta(P^2,\sgn P;-1)=({\rm
ord}\,P)\Res_1\eta(P;0).
  \mylabel{Grev-2.34}
\end{equation}

2. To prove the converse we consider a classical pseudodifferential idempotent,
$B$, on a compact manifold, $M$. $B$ is similar to a self--adjoint idempotent
and it is not difficult to see that the similarity can be effected
through a pseudodifferential operator. Since the residue
is a trace, similar operators have the same residue. Hence we may assume
$B$ to be an orthogonal projection. The assertion will follow from
\myref{Grev-2.34} if we can show that there exists an invertible self--adjoint
classical pseudodifferential operator, $P$, of order $1$ with
$$B=\frac 12(\sgn P+I).$$
We choose a first order self--adjoint classical pseudodifferential
operator, $Q$, with scalar principal symbol $\sigma_Q(\xi)=|\xi|$.
Furthermore, we may choose $Q$ to be positive. Then put
$$\tilde P:=BQB-(I-B)Q(I-B).$$
$\tilde P$ is elliptic and commutes with $B$. To make it invertible we put
$$Px:=\left\{\begin{array}{rl}
    \tilde Px,&\quad x\in \ker \tilde P^\perp,\\
    x,&\quad x\in\ker\tilde P\cap\Im B,\\
   -x,&\quad x\in\ker\tilde P\cap\ker B.
   \end{array}\right.$$
By construction we have $B=\frac 12(\sgn P+I)$ and hence we reach
the conclusion.\endproof

\comment{\noindent
2. To prove the converse we consider a classical pseudodifferential idempotent,
$B$, on a compact manifold, $M$. $B$ is similar to a self--adjoint idempotent
and it is not difficult to see that the similarity can be effected
through a pseudodifferential operator. Since the residue
is a trace, similar operators have the same residue. Hence we may assume
$B$ to be an orthogonal projection. We put $\tau:=2P-I$ and let
$\sigma_\tau\in\cinf{S^*M}$ be the principal symbol of $\tau$.
We can choose an invertible first order self--adjoint pseudodifferential
operator, $Q$, with principal symbol $\sigma_\tau$. Then we put
$$P:= \frac 12 ( Q|Q|^{-1}+I).$$
This is an operator of order $0$ and $P-B$ is of order $-1$.
Then one shows that there exists a pseudodifferential
projection $P_1$, a smoothing operator $R$, and a pseudodifferential
operator $K$, $\|K\|<1$,
such that $P+R$ is a projection and
$$B=P+R+K.$$
Since $\|K\|<1$ the projections $B$ and $P+R$ are similar and since
\res is a trace which vanishes on smoothing operators we find
$$\res B=\res(P+R)=\res P.$$
Since $\res I=0$ we end up with
\begin{epeqnarray}
    \res P&=&\frac 12 \res(Q|Q|^{-1})\\
         &=& \frac 12 \Res_1\eta(Q;0)=0.
\end{epeqnarray}         }
  
\vspace{0.5em}\noindent
We emphasize, however, that neither for index theorems \cite{BS} nor for
the gluing law to be proved below the regularity at $0$ of the
$\eta$--function is essential; the
definition \myref{G4-1.12} is perfectly sufficient.

If one wants to widen the class of operators which admit reasonable
$\eta$--invariants then it is most natural to consider elliptic
boundary value problems. As illustrated by the gluing question,
one may also expect further insight in the compact case. The first
work in this direction seems to be \cite{GilkeySmith} which deals
with local boundary conditions leading to (mildly) nonself--adjoint
operators which do, however, admit reasonable $\eta$--invariants. This
was used by Singer \cite{Singer} who showed (among other things) that
the difference of $\eta$--invariants associated to two natural
boundary value problems of this kind is an interesting spectral invariant of
the boundary, at least asymptotically. More precisely, let $M$
be an odd dimensional Riemannian spin manifold with spinor bundle
$S(M)$ and assume again that the metric is a product near $N$
(this assumption will be kept from now on). Thus, a neighborhood
of $N$ in $M$ is isometric to the cylinder $N_R=[0,R)\times N$, 
for some $R>0$, therefore we will write $M_R$ to make the dependence on
$R$ more transparent.\revmark Then we have again a representation of type
\myref{G1-1.2} for the Dirac operator, $D^{M_R}$, on $S(M_R)$ where
$A=D^N$ becomes the Dirac operator on $S(N)=S(M_R)|N$. Under
$\gamma$, $S(N)$ splits into $S^+(N)\oplus S^-(N)$ with projections
$Q_\pm:L^2(S(N))\longrightarrow L^2(S^\pm(N))$. Then
$D_\pm^{M_R}:= (D^{M_R},Q_\pm)$ are well--posed boundary value 
problems to which
the analysis of \cite{GilkeySmith} applies, and Singer proves
that by stretching $N_R$ the difference of $\eta$--invariants localizes
i.e.
\begin{equation}
  \lim_{R\to\infty} (\eta(D_+^{M_R})-\eta(D_-^{M_R}))=\frac{1}{4\pi i}
      \log\det(D^N)^2.\mylabel{G4-1.22}
\end{equation}
\revmark
Singers investigation was motivated by Witten's identification of the covariant
anom\-a\-ly with the so--called adiabatic limit of an $\eta$--invariant \cite{Witten}
but his work, in turn, stimulated greatly the interest in $\eta$--invariants
for manifolds with boundary.

Douglas and Wojciechowski \cite{Dougwoj} then studied systematically the
properties of $\eta$--invariants for generalized Dirac operators
on odd--dimensional manifolds with boundary. They assumed 
\myref{G1-1.2} with the additional hypothesis
\begin{equation}
    \ker A=0,\mylabel{G4-1.23}
\end{equation}
and chose the boundary condition \myref{G4-1.3a}; in this situation,
they established Lemmas \plref{S4-1.1} and \plref{S4-1.2}, and
for suitable families of such operators they proved
\myref{G4-1.18} for $k=0$. Moreover, they showed that stretching
the cylinder $N_R$ produces an "adiabatic limit" in the sense that
\begin{equation}
   \lim_{R\to\infty} \eta(D_R)=:\eta_\infty
   \mylabel{G4-1.24}
\end{equation}
exists.
Then the challenge was to identify $\eta_\infty$ and to extend the
results to $\ker A\neq 0$. In this case, there is considerable freedom
of choice for the "supplementary" boundary condition (\plref{G4-1.4a},b),
and its variation ought to be allowed, too, in a suitable generalization
of \myref{G4-1.18}. Note that the analysis of Lemma
\plref{S4-1.2} does not apply to this situation right away since the
operators under consideration do not have constant domain, so one
has to search for a suitable transformation of the family.
This was done by Lesch and Wojciechowski \cite{LeschW}.
Since their method also served as a basic motivation for this paper, we will
present a suitable version of their argument.
Theorem \plref{S6-2.5} below generalizes considerably the original
construction and is the main analytic tool of our present work.

The result of \cite{LeschW}
was obtained independently by M\"uller \cite{Muller}. In addition,
M\"uller presented a thorough analysis of the operators
$D_\sigma$ in the general case. In particular, he showed that
$\eta_\infty$ exists and can be interpreted as the 
suitably defined $\eta$--invariant for an operator on the
manifold $\tilde M:=M\cup N_\infty$. Moreover, he proved that
\begin{equation}
   \eta_\infty=\eta(D_{\sigma_1})
   \mylabel{G4-1.24.5}
\end{equation}
for a suitable $\sigma_1$, obtained from scattering theory on
$\tilde M$. He also obtained the regularity of the $\eta$--function
of $D_\sigma$ if $D$ is assumed to be of Dirac type.

In the context of Melrose's "b--calculus",
Hassell, Mazzeo and Melrose \cite{MazMel,MazMelHas} define
an $\eta$--invariant on manifolds with boundary, and they
prove a gluing law in this situation. This $\eta$--invariant
coincides again with $\eta_\infty$.

\myref{G4-1.24.5} can be taken as the starting point to prove the
gluing law for $\eta$--invariants as done by M\"uller
\cite{Muller1}
\revmark\ and Wojciechowski \cite{Woj1, Woj2}.
Bunke \cite{Bunke} gave a complete proof of the gluing law based on
cutting the manifold in question thrice and reassembling
the pieces into a cylinder (carrying both boundary conditions)
and a compact manifold where one can do essentially
only "interior" analysis, in view of the finite propagation speed
enjoyed by all $D_\sigma$. This reduces the analysis
to the explicit computation on the cylinder
carried out in \cite{LeschW}. Bunke's result is, at least theoretically,
more precise than ours since he gives a formula for the unknown
integer in \myref{G4-1.20}. This is possible since his deformation
induces a relatively compact perturbation. By contrast, our construction
is more direct and more general but less rigid with regard to
compactness.

Bunke's argument, in turn, was generalized and simplified in a substantial
paper
by Dai and Freed \cite{DF}; they interpreted the invariant
\myref{G4-1.21} as a section of the determinant line if one
considers families of operators $D_\sigma$ fibered over a compact
Riemannian manifold. This allows a natural interpretation of Witten's
anomaly formula, and also illustrates nicely the philosophy
developed in Singer's paper \cite{Singer}.

Our proof of the gluing law (Theorem \plref{S7-2S.3} below) arises
as a byproduct of an extension of the variation formula
to a wider class of boundary conditions, thus furnishing
a proof of a rather different nature than those 
described before.

\section{Expansion theorems and the gluing law}\label{Sec2}

Our approach to the proof of the gluing law was originally
inspired by Vishik's proof of the Cheeger--M\"uller Theorem
\cite{Vishik}. Working out the details we discovered, however,
that we were lead to a very natural generalization
of the approach in \cite{LeschW}, designed to determine
the variation of $\eta(D_\sigma)$ under a change of $\sigma$.

At any rate, the analysis we are going to present
deals with operators of type \myref{G1-1.2} but with more general
boundary conditions than \G4-1.3. We will now explain how this
class arises naturally from the gluing problem, define it in general,
and outline the proof of the gluing law. Most details are deferred to
Sec. \mls\plref{Sec3}

Let now  $M$ be a {\it compact} Riemannian manifold, $\dim M=m$, and
let
\begin{equation}
 D_0:\cinfz{S}\longrightarrow \cinfz{S}
 \mylabel{G5-2.1}
\end{equation}
be a first order symmetric elliptic differential operator on the
hermitian vector bundle $S\to M$. The main examples
are, of course,
Dirac operators associated to a Dirac bundle $(S,\nabla)$, but we will
work in a more general context, allowing for example
Dirac operators with potential.

Let $N\subset M$ be a compact
hypersurface. We assume that $N$ has a tubular
neighborhood $U$ isometric to $(-1,1)\times N$ and such that the
hermitian structure of $S$ is a product, too. Moreover, we assume
that on $U$ the operator $D_0$ has the form
\begin{equation}
D_0=\gamma(\frac{\pl}{\pl x}+A),
  \mylabel{G5-2.2}
\end{equation}
where $\gamma\in\cinf{\End(S_N)}$ is a unitary bundle automorphism and
$A$ is a first order self--adjoint elliptic differential operator on
$S_N:=S|N$.
If $D_0$ is a compatible Dirac operator,  then $\gamma$ is Clifford
multiplication by the inward normal vector and $A$ is (essentially)
a Dirac operator on $N$. We assume, furthermore, that $\gamma$
and $A$ satisfy \myref{G1-1.3}.

Let $D$ be the restriction of $D_0$ to $\cinfz{S|M\setminus N}$. This
operator is no longer essentially self--adjoint;
in order to obtain self--adjoint extensions
one has to impose boundary conditions.
The natural boundary condition inherited from $M$
is the {\it continuous transmission} boundary
condition. Interpreting sections of $S$ with support in $U$ as functions
$[-1,1]\to L^2(S_N)$ in the obvious way, this
boundary condition reads
\begin{equation}
    f(0-)=f(0+).
    \mylabel{G1-1.4}
\end{equation}
It is fairly clear that the resulting self--adjoint operator is unitarily
equivalent to the closure of $D$ in $L^2(S)$. On the other hand,
$D$ lives naturally on
\begin{equation}
   M^{\rm cut}:=(M\setminus U)\cup_{\pl(M\setminus U)}
     \big( (-1,0]\times N \cup [0,1)\times N\big )
     \mylabel{G1-1.5}
\end{equation}
obtained by cutting $M$ along $N$ (we adopt here the notation
from \cite[p. \mls 5164 and Sec. \mls 4]{DF}). Thus,
$\mcut$ is obtained from $M$ by artificially introducing two
copies of $N$ as boundary.

On $\mcut$ we can introduce spectral boundary conditions as in 
Sec. \mls\plref{Sec1}
The natural interpolation between the continuous transmission
and the Atiyah--Patodi--Singer boundary condition is furnished
by the boundary conditions
\alpheqn
\begin{eqnarray}&\begin{array}{rcl}\DST
     \cos \theta\, P_{>0}(A) f(0+)&=& \DST\sin \theta\, P_{>0}(A)f(0-),\\[0.5em]
     \DST\sin\theta\, P_{<0}(A)f(0+)&=&\DST \cos \theta\, P_{<0}(A)f(0-),
     \end{array}&\mylabel{G1-1.6}\\[0.5em]
     &P_0(A) f(0+)= P_0(A) f(0-), \mylabel{G5-2.5b}
\end{eqnarray}
\reseteqn
where $|\theta|<\pi/2$.

To render this more transparent, we employ the isomorphism
(with $H:=L^2(S_N)$)
\alpheqn
\begin{equation}
   \Phi:L^2(S|U) \simeq L^2([-1,1],H)\longrightarrow
      L^2([0,1],H\oplus H),
\end{equation}
which sends $f\in L^2([-1,1],H)$ to $\Phi f$,
\begin{equation}
  \Phi f(x)= f(x)\oplus f(-x), \quad x\in [0,1].
\end{equation}
\reseteqn 
It is easy to see that, under $\Phi$, $D$ is transformed to
\begin{equation}
      \tilde D:=
      \mat{\gamma}{0}{0}{-\gamma}\left(\frac{\pl}{\pl x}+\mat{A}{0}{0}{-A}\right)
         =: \tilde\gamma\Big(\frac{\pl}{\pl x} +\tilde A\Big),
         \mylabel{G1-1.8}
\end{equation}
and the boundary condition to
\alpheqn
\begin{equation}
  \cos\theta\, P_{>0}(\tilde A) u(0)= \sin \theta\; \tau P_{<0}(\tilde A)u(0),
\end{equation}  
where
\begin{equation}
    \tau=\mat{0}{1}{1}{0}\otimes I_H,
    \mylabel{G1-1.10}
\end{equation}
supplemented on $\ker\tilde A$ by
\begin{equation}
   P_{\sigma} u(0)=0,
\end{equation}
with
\begin{equation}
   \sigma:= -\mat{0}{P_0(A)}{P_0(A)}{0}.
   \mylabel{Grev-3.8d}
\end{equation}   
\reseteqn
Note that
\begin{equation}
\tau\tilde\gamma+\tilde\gamma\tau=0=\tau\tilde A+\tilde A\tau,\quad
  \tau^2=1,\quad \tau=\tau^*.\mylabel{G5-2.9}
\end{equation}   

Next we observe that this boundary condition can be written as
\alpheqn
\begin{equation}
   \tilde P(\theta) u(0)=0,
   \mylabel{G6-2.10a}
\end{equation}
if we introduce the projection
\begin{equation}
 \tilde P(\theta):= \cos^2\theta P_{>0}(\tilde A)+\sin^2\theta P_{<0}(\tilde A)
    -\frac 12(\sin 2\theta) \tau (P_{>0}(\tilde A)+P_{<0}(\tilde A))+ P_{\sigma}.
    \mylabel{G6-2.10b}
\end{equation}
\reseteqn

\noindent
It is useful to note the following properties of this family of
projections, all of which are easily verified.

First, we see that
\begin{equation}
   \tilde\gamma \tilde P(\theta)=(I-\tilde P(\theta)) \tilde\gamma,
   \mylabel{G5-2.11}
\end{equation}
and that $\tilde P(\theta)$ commutes with $\tilde A^2$,
\begin{equation}
   [\tilde P(\theta),\tilde A^2]=0.
\end{equation}
We do not have commutativity with $\tilde A$, however. Instead
we find
\begin{equation}
   \tilde P(\theta) \tilde A \tilde P(\theta)= \cos 2\theta |\tilde A|
     \tilde P(\theta).\mylabel{G5-2.13}
\end{equation}

Remembering the argument of Lesch and Wojciechowski \cite{LeschW}
we are lead to ask for a natural ''parametrization'' of the family
$(\tilde P(\theta))_{|\theta|<\pi/2}$. It is easy to verify that
with
\begin{equation}
   U(\theta):=\left(\cos \theta (P_{>0}(\tilde A)+P_{<0}(\tilde A))
         +\sin\theta(P_{>0}(\tilde A)-P_{<0}(\tilde A))\tau\right)
         \oplus I_{\ker\tilde A}
         \mylabel{G5-2.14}
\end{equation}
and
\begin{equation}
  \sgn \tilde A:= P_{>0}(\tilde A)- P_{<0}(\tilde A)
\end{equation}
we have
\begin{eqnarray}
    \tilde P(\theta)&=& U(\theta) \tilde P(0) U(\theta)^*,\mylabel{G5-2.16}\\
    U(\theta)&=& e^{(\sgn \tilde A \tau) \theta}.\mylabel{G5-2.17}
\end{eqnarray}

Thus we obtain a family of generalized Atiyah--Patodi--Singer
boundary conditions, and the gluing law becomes just the variational
formula for this class of operators in the sense of Sec. \mls\plref{Sec1}

In fact, we will generalize the situation further. {\it Thus from
now on we consider the following setting.}

$M$ is a Riemannian manifold  of dimension $m$, $S\to M$ is a 
smooth hermitian vector bundle over $M$, and $D$ is a first order
symmetric elliptic differential operator on $\cinfz{S}$. We
assume that $M$ can be decomposed as 
\begin{equation}
   M=U\cup M_1,
   \mylabel{G5-2.18}
\end{equation}
where $M_1$ is a compact manifold with boundary $N=\partial M_1=\partial U$
and $U$ is open. Moreover, we assume an isometry of Hilbert spaces,
\begin{equation}
   \Phi: L^2(S|U)\longrightarrow L^2([0,1],H),\mylabel{G5-2.19}
\end{equation}
where $S_N$ is a smooth hermitian bundle over $N$ and $H=L^2(S_N)$ as before.
This isometry
maps smooth sections to smooth sections in the sense
that
\begin{equation}
 \Phi(\cinf{S|U}\cap L^2(S|U))\subset \cinf{(0,1),\cinf{S_N}}\cap
      L^2([0,1],H).
\end{equation}
Thus we can transform $D$ on $U$, and we require that
\begin{equation}
    \Phi D \Phi^* =\gamma(\pl_x +A)=:\tilde D,
\end{equation}                          
with $A$ a symmetric elliptic operator of first order on $S_N$ which
we identify with its self--adjoint closure,
and $\gamma$ a bounded operator on $H$. We assume,
moreover, that $\gamma$ and $A$ satisfy the relations
\myref{G1-1.3} and \myref{G4-1.7}.

Finally, we require that for $\phi\in\cinfz{-1,1}$ there is
$\psi_\phi\in\cinf{M}$ such that $\psi_\phi=0$ in a neighborhood
of $\partial M_1$, and
\alpheqn
\begin{equation}
   \Phi(\psi_\phi u)=\phi \Phi u,\quad u\in L^2(S);
\end{equation}
and
\begin{equation}
   \phi=1 \enspace\mbox{\rm near}\enspace 0 \enspace\mbox{\rm implies}\enspace
    1-\psi_\phi \in\cinfz{M}.\mylabel{G5-2.22b}
\end{equation}\reseteqn
As usual, we extend $\tilde D$ to $L^2(\R_+,H)=:\ch$
to obtain the model operator. To define a family of boundary
conditions we proceed as in the above analysis of the cutting
problem: we consider a family ${P(\theta)}_{|\theta|<\pi/2}$
of orthogonal projections with the following properties.

\begin{numabsatz}
\myitem{ $\gamma P(\theta) = (I-P(\theta)) \gamma$;\mylabel{G5-1.31a}}
\myitem{ $[P(\theta),A^2]=0$;\mylabel{G5-1.31b}}
\myitem{ $A(\theta):= P(\theta)AP(\theta)=a(\theta) |A|P(\theta)$ for some\\
  $a\in\cinf{-\pi/2,\pi/2}$ with $a>-1$.\mylabel{G5-2.25}}
\end{numabsatz}  
These projections are again assumed to be
conjugate to $P(0)$ under a family of unitaries, $U(\theta)$,
\begin{equation}
  P(\theta)=U(\theta) P(0) U(\theta)^*.\mylabel{G5-2.26}
\end{equation}
We assume, moreover, a representation
\begin{equation}
   U(\theta)= e^{i T(\theta)},\mylabel{G5-2.27}
\end{equation}
with $T(\theta)$ bounded and self--adjoint in $H$, smooth in
$(-\pi/2,\pi/2)$, and such that
\alpheqn
\begin{eqnarray}
    &&[\gamma,T(\theta)]=0,\mylabel{G5-2.28a}\\[0.5em]
   &&A T(\theta)+T(\theta) A =0.\mylabel{G5-2.28}
\end{eqnarray}
\reseteqn

With these data we define boundary conditions for $D$ and $\tilde D$ via
\revmark
\alpheqn   
\begin{eqnarray}
  \tilde \cd_\theta&:=&\Big\{ u\in C(\R_+,H)\cap \ch\,\Big|\,
     u\in\cd(\tilde D^*),  P(\theta)u(0)=0\Big\},\\[0.5em]
         \cd_\theta&:=&\left\{ u\in L^2(S)\,\Bigg|\,
         \begin{array}{c} u\in \cd(D^*), \Phi (\psi_\phi u)\in\tilde \cd_\theta
         \quad\mbox{for some}\\[0.5em]
         \phi\in\cinfz{-1,1}\quad\mbox{with}\quad \phi=1\quad\mbox{near}\quad 0
         \end{array}\right\},
\end{eqnarray}   
\reseteqn
and
\begin{equation}
   D_\theta:=D|\cd_\theta,\quad \tilde D_\theta:=\tilde D|\tilde\cd_\theta.
\end{equation}

A good part of the subsequent analysis rests on these assumptions. For
the asymptotic expansions to exist it is convenient to require
in addition that
\begin{equation}
\begin{array}{l}
\mbox{$P(\theta), T(\theta)$ are classical pseudodifferential operators}\\
\mbox{of order zero on $N$, for $|\theta|<\pi/2$.}
\end{array}\mylabel{G6-2.31}
\end{equation}
This assumption is clearly satisfied in the gluing case
{\rm (\ref{G6-2.10a},b)}.

We will refer to the family $(D_\theta)_{|\theta|<\pi/2}$ with the properties
listed above as a {\it deformation of Atiyah--Patodi--Singer} (APS) type.
Then we have seen that cutting along a compact hypersurface leads naturally
to such a family. In this case, we do have a bit more structure since, in
\myref{G5-2.25}, we have $a(\theta)=\cos 2\theta$, in view of
\myref{G5-2.13}, and we have the additional symmetry, $\tau$, with the
properties \myref{G5-2.9}.

We note that a single projection, $P$, with the properties
\myref{G5-1.31a}, \myref{G5-1.31b}, \myref{G5-2.25}  defines a self--adjoint
extension of $D$, $D_P$, to which the analysis of Sec. \mls\sectone\            
applies. This we call a {\it generalized APS operator} since, clearly, 
$P=P_{>0}(A)+P_\sigma$ falls in this class.

We proceed to the spectral analysis of $D_\theta$, the proofs being given
in Sec. \mls\plref{Sec3}
\begin{prop}{S5-2.1} The operators $D_\theta$ and $\tilde D_\theta$ are
essentially self--adjoint.
\end{prop}

\bigskip\noindent
We will identify $D_\theta$ and $\tilde D_\theta$ with their respective
closures in the sequel.

\begin{prop}{S5-2.2} $D_\theta$ satisfies \myref{G4-1.7} i.e.
$$(D_\theta+i)^{-1}\in C_p(L^2(S))\quad\mbox{for every}\quad p>m.$$
\end{prop}

\bigskip\noindent
We want to apply Lemma \plref{S4-1.3.5} to the family $(D_\theta)_{|\theta|<\pi/2}$
which requires that we first apply a transformation to
satisfy (\ref{G4-1.13a},e). This we do as in \cite{LeschW}, and this is
the motivation for the assumptions \myref{G5-2.26}, \myref{G5-2.27}, and
(\ref{G5-2.28a},b).

Thus we choose $\phi\in\cinfz{-1,1}$ with $\phi=1$ near $0$ and
introduce the unitary transformation
\begin{equation}\begin{array}{c}\DST
  \Psi_\theta:L^2([0,1],H)\longrightarrow L^2([0,1],H),\\[0.5em]
  \DST\Psi_\theta u(x):=e^{i \phi(x) T(\theta)} (u(x)).
  \end{array}\mylabel{G5-2.31}
\end{equation}
Then $P(0)u(0)=0$ implies $P(\theta) \Psi_\theta u(0)=0$, in view of
\myref{G5-2.26}. Hence, extending $\Psi_\theta$ to $L^2([0,1],H)\oplus
L^2(S|M_1)$ as the identity on $L^2(S|M_1)$ and similarly $\Phi$ in
\myref{G5-2.19}, we obtain an isometry
 $$\Phi_\theta:=\Phi^* \Psi_\theta \Phi$$
of $L^2(S)$ mapping $\cd_0$ to $\cd_\theta$. Consequently, the family
\begin{equation}
\check D_\theta:=\Phi_\theta^*D_\theta \Phi_\theta
  \mylabel{G5-2.32}
\end{equation} 
has constant domain, $\cd_0$, and the same spectral invariants
as $D_\theta$. It is easy to see that $(\check D_\theta)_{|\theta|<\pi/2}$ satisfies
(\ref{G4-1.13a},e). It remains to establish the asymptotic expansions
\myref{G4-1.11}, with   $\check D_\theta, \frac{d}{d\theta}\check D_\theta$ in place
of $B$.

Our expansion results will be expressed in terms of the Mellin
transform of a certain meromorphic
function, $F_a$, which we have to introduce first.

\begin{lemma}{S6-3.2} Consider for $a\in(-1,1]$ and $x>0$ the
function
$$F_a(x):=x \int_0^\infty \erfc(z) e^{-2axz-x^2} dz,$$
where
$$\erfc(z)=\frac{2}{\sqrt{\pi}}\int_z^\infty e^{-u^2} du$$
denotes the complementary error function.

Then the Mellin transform of $F_a$ is, for $0<|a|<1$,
\begin{eqnarray}
  \cm F_a(w)&=&\frac{1}{4a}\Big[(1-(1-a^2)^{-w/2})\Gamma(w/2)\nonumber\\
      &&+\frac{2}{\sqrt{\pi}}(1-a^2)^{-w/2}
         \int_0^a (1-t^2)^{(w/2)-1} dt\; \Gamma((w+1)/2)\Big],
          \mylabel{G6-L32.1}
\end{eqnarray}
whereas
\begin{equation}
    {\cm}F_0(w)=\frac{1}{2\sqrt{\pi}}\Gamma((w+1)/2),
            \mylabel{G6-L32.2}
\end{equation}
and
\begin{equation}
           {\cm}F_1(w)=\frac 14\Big[\Gamma(w/2)-\frac{2}{w\sqrt{\pi}}
            \Gamma((w+1)/2)\Big].
                     \mylabel{G6-L32.3}
\end{equation}
Hence ${\cm}F_a(w)$ is meromorphic in $\C$ with simple poles at the
points $-k, k\in\Z_+$. For $|a|<1,$ the residues are
\begin{equation}
    \begin{array}{rcl}
     \Res_1 {\cm}F_a(-2l)&=&\DST \frac{(-1)^l}{l! 2a}(1- (1-a^2)^l), \quad
           l\in\Z_+,\\[1em]
     \Res_1 {\cm}F_a(-2l-1)&=&\DST \frac{(-1)^l}{l!\sqrt{\pi}a}
          (1-a^2)^{l+1/2}\int_0^a (1-t^2)^{-l-3/2} dt, \quad
          l\in\Z_+.
     \end{array}
     \mylabel{G6-L32.4}
\end{equation}

For $a=0,1$ one has to take the corresponding limit in \myref{G6-L32.4}.
More precisely,
\begin{equation}
     \widearray
     \begin{array}{rcl}
       \Res_1 {\cm}F_0(-2l)&=&\DST 0,\quad l\in \Z_+,\\[0.5em]
       \Res_1 {\cm}F_0(-2l-1)&=&\DST \frac{(-1)^l}{\sqrt{\pi}l!},\quad l\in \Z_+,\\[1em]
       \Res_1 {\cm}F_1(-2l)&=&\DST \casetwo{0}{l=0,}{\DST \frac{(-1)^l}{l!2}}{l\in\N,}\\[1.5em]
       \Res_1 {\cm}F_1(-2l-1)&=&\DST \frac{(-1)^l}{l!\sqrt{\pi}(2l+1)},\quad l\in \Z_+.
     \end{array}
     \mylabel{G6-L32.5}
\end{equation}
\end{lemma}

\vspace{1em}\noindent
Now we present our first expansion result.

\begin{theorem}{S1-2.7}
Assume that \myref{G5-2.18} through \myref{G6-2.31} hold.
For $l=0,1$ we have an
asymptotic expansion of the form
\begin{eqnarray}
     \tr_{L^2(S)}\big[\dth^l\edth\big]&\sim_{t\to 0+}&
      \sum_{j=0}^\infty a_j(\theta,l)\;t^{j-m/2}
      +\sum_{j=0}^\infty b_j(\theta,l)\;t^{j/2}\,\log t\nonumber\\
      &&+\sum_{j=0}^\infty c_j(\theta,l)\;t^{(j-n-l)/2}
        +\sum_{j=0}^\infty d_j(\theta,l)\;t^{j/2}.
     \mylabel{G1-2.26}
\end{eqnarray}
Here, the coefficients $a_j$ are integrals of local densities
on the metric double, $\tilde M$, of $M$, $b_j$ and
$c_j$ are integrals of local densities on $N$, and
$d_j$ are nonlocal invariants of $N$; they are given
explicitly in the formulas
\myref{G2-2.29}, \myref{G6-3.14a}, \myref{G6-3.14b},
\myref{G6-3.24}, and \myref{G6-3.26} below.

For $l=0$, the leading term is
\begin{equation}
    a_0(\theta,0)=\Gamma(m/2+1) \vol (T_1^*M),
    \mylabel{G6-2.35}
\end{equation}
where $T_1^*M=\{\xi\in T^*M\,|\, \sigma_{D_\theta^2}(\xi)\le 1\}.$

The logarithmic terms vanish if $l=0$ and $m$ is odd.
If $l=0$ and $m$ is even then $b_{2j}(\theta,0)=0$. However,
the logarithmic terms are present in general.

For $l=1$, the expansion \myref{G1-2.26} implies that
$\eta(D_\theta;s)$ has a meromorphic extension to $\C$
with at most double poles.
$0$ is a simple pole and
for the residue at $0$ we find
\begin{equation}
   \Res_1 \eta(D_\theta;0)= \frac{2}{\sqrt{\pi}} a_{n/2}(\theta,1)
      +\Big(\frac{2a(\theta)}{\sqrt{\pi}} {\cm}F_{a(\theta)}(1)-\frac 12\Big)
        \res (\gamma (\sgn A) P(\theta)).
\mylabel{G6-2.36}
\end{equation}
\end{theorem}

\vspace{1em}\noindent
For the APS boundary condition, this result has been obtained
by Grubb and Seeley \cite{GrubbSeeley}. 
\comment{Our approach differs from theirs
by using the explicit heat kernel \myref{G5-3.1}; this method seems to
provide explicit formulas for the coefficients in
\myref{G1-2.26} more directly.}

By contrast, our approach is simply based on the spectral theorem and
the explicit formula \myref{G5-3.1}. Nevertheless,  we can handle
boundary conditions which are significantly different from the 
APS condition.

 The expansion result for the APS condition is sketched
in M\"uller \cite[Lemma 1.17]{Muller} overlooking,
however, the coefficients which are not local in $\tilde M$ in the case
$l=0$. In the case $l=1$ and for APS boundary conditions,
these nonlocal terms are actually not present.

To explain this let for the moment $\tilde D_\sigma$ be the operator with
APS boundary condition. Then a simple symmetry argument shows that for any
cut--off function $\phi\in\cinfz{\R}$ 
\begin{equation}
\tr_{L^2(S)}\big[\phi  \tilde D_\sigma e^{-t \tilde D_\sigma^2}\big]=0,\quad t>0,
\mylabel{Grev-3.42}
\end{equation}
and hence $b_j(\sigma,1)=c_j(\sigma,1)=d_j(\sigma,1)=0$
(cf. \cite[Lemma 5.2.4]{Lesch2}). For general $P(\theta)$ we cannot
expect \myref{Grev-3.42} to hold.

\bigskip
In the next step, we evaluate the formula for the variation of
the $\ovl{\eta}$--invariant in Lemma \plref{S4-1.3.5}, via the asymptotic
expansion of $\tr\big((\frac{d}{d\theta} \check\dth) e^{-t\check\dth^2}\big)$.

\begin{theorem}{S6-2.5}
Under the assumptions of Theorem \plref{S1-2.7} we have the following
variation formulas:
\begin{eqnarray}
   \frac{d}{d\theta} \Res_1\eta(D_\theta;0)&=&\frac{1}{\sqrt{\pi}}
            \res(\gamma i T'(\theta)),\\
   \frac{d}{d\theta} \ovl{\eta}(D_\theta)&=&
        \frac{1}{2\pi} a_{00}(A,\gamma i T'(\theta))\nonumber\\
        &&+\Big(\frac{2a(\theta)}{\sqrt{\pi}} {\cm}F_{a(\theta)}(1)-\frac 12\Big)
        \res(\gamma i T'(\theta) (\sgn A) P(\theta)).
\end{eqnarray}
\end{theorem}

\bigskip

\begin{cor}{S4-1.6} In the situation of Theorem \plref{S6-2.5},
assume in addition that
\begin{equation}
T'(\theta)P(\theta)=(I-P(\theta)) T'(\theta).
\mylabel{G7-2S.0}
\end{equation}
Then
$$\res(\gamma i T'(\theta) (\sgn A) P(\theta))
  = \frac{a(\theta)}{2} \res(\gamma i T'(\theta)).$$
In particular, if $\res(\gamma i T'(\theta))=0$ then
$\Res_1\eta(D_\theta;0)$ is independent of $\theta$ and
$$\frac{d}{d\theta} \ovl{\eta}(D_\theta)=\frac{1}{2\pi}
   a_{0,0}(A,\gamma i T'(\theta)).$$
\end{cor}

\proof We use \myref{G7-2S.0}, \myref{G5-1.31a}, \myref{G5-2.25},
and the trace property
of the noncommutative residue to compute
\begin{eqnarray*}
      \res(\gamma i T'(\theta) (\sgn A) P(\theta))&=&
      \res(\gamma i T'(\theta)P(\theta) (\sgn A) P(\theta))\\
      &=& a(\theta) \res(\gamma i T'(\theta)P(\theta)).
\end{eqnarray*}
Here we have used that $\res$ vanishes on smoothing operators. Furthermore,
in view of \myref{G5-2.28a},
\begin{eqnarray*}
      \res(\gamma i T'(\theta) P(\theta))&=&
      \res(\gamma i(I-P(\theta)) T'(\theta))\\
      &=&\res(i(I-P(\theta)) T'(\theta)\gamma)\\
      &=&\res(i(I-P(\theta)) \gamma T'(\theta))\\
      &=& \res(\gamma i T'(\theta)(I-P(\theta))),
\end{eqnarray*}
and we reach the conclusion.\endproof

\noindent
Next we introduce a special class of deformations of APS type which is
still slightly more general than the gluing situation
\myref{G1-1.6}--\myref{G5-2.17}:

We consider again the framework \myref{G5-2.18}--\myref{G5-2.22b}.
Furthermore, let
$\tau:\cinf{S_N}\to \cinf{S_N}$
be a unitary classical pseudodifferential operator
satisfying (cf. \myref{G5-2.9})
\begin{equation}
    \tau\gamma+\gamma\tau=0=\tau A+A\tau,\quad \tau^2=I, \quad
    \tau=\tau^*.
   \mylabel{G7-2S.1}
\end{equation}
We abbreviate
\begin{equation}
       K^\pm:=(\ker A)\cap \ker(\gamma\mp i).
\end{equation}
The relations \myref{G7-2S.1} immediately imply
\begin{equation}
    \dim K^+=\dim K^-.\mylabel{G7-2S.2}
\end{equation}

However, the presence of $\tau$ is not really necessary for this
equality.  \myref{G7-2S.2} follows already from
\myref{G5-2.18}--\myref{G5-2.22b}.
If $D$ is a Dirac
operator, this is the well--known cobordism theorem for
Dirac operators \cite[Chapter XVII]{Palais}. For general $D$, this is
due to the second named author \cite[Theorem 6.2]{Lesch1},
\cite[Chapter IV]{Lesch2}. It was also proved independently by
W. M\"uller \cite[Prop. \mls 4.26]{Muller}.\label{cobordism}

In view of \myref{G7-2S.2} we can choose an isometry
\begin{equation}
    U:K^+\longrightarrow K^-
    \mylabel{G7-2S.3}
\end{equation}
and put
\begin{equation}
    \sigma=\mat{0}{U^*}{U}{0}:\ker A\longrightarrow \ker A.
    \mylabel{G7-2S.4}
\end{equation}

With these data we can introduce the projection (cf. \myref{G6-2.10b})
\begin{equation}
 P(\theta):= \cos^2\theta P_{>0}(A)+\sin^2\theta P_{<0}(A)
    -\frac 12(\sin 2\theta) \tau (P_{>0}(A)+P_{<0}(A))+ P_{\sigma}
    \mylabel{G7-2S.5}
\end{equation}
and the unitary family (cf. \myref{G5-2.14})
\begin{equation}
   U(\theta)=\big(\cos \theta (P_{>0}(A)+P_{<0}(A))
         +\sin\theta(\sgn A)\tau\big)
         \oplus I_{\ker A}=e^{(\sgn A) \tau \theta}.\mylabel{G7-2S.6}
\end{equation}
One immediately checks the relations \myref{G5-2.11}--\myref{G5-2.13},
\myref{G5-2.16}, hence we are lead to a deformation of APS type. We
denote the corresponding family of operators by $D_{\theta, \sigma}$,
indicating explicitly the dependence on the choice of $\sigma$.
If we fix $\theta$ and consider a one parameter family of
reflections, $\sigma_u$, we obtain another deformation of APS type.
In this way we recover the main result of Lesch and Wojciechowski
\cite{LeschW} as a special case of our present work:

\begin{prop}{S7-2S.1} {\rm(cf. \cite{LeschW,Muller,DF}) }
Let $\cos\theta\not=0$ and $U_u:K^+\to K^-$
be a smooth family of unitary operators. Put
$$\sigma_u:=\mat{0}{U_u^*}{U_u}{0}.$$
Then $(D_{\theta,\sigma_u})_u$ is a deformation of {\rm APS} type,
$\Res_1 \eta(D_{\theta,\sigma_u};0)$ is independent of $u$ and
$$\frac{d}{du} \ovl{\eta}(D_{\theta,\sigma_u})=
  \frac{1}{2\pi i} \tr_{K^+}[U_u^{-1}\frac{d}{du} U_u].$$
\end{prop}
\proof We put
$$P_u(\theta):=\cos^2\theta P_{>0}(A)+\sin^2\theta P_{<0}(A)
    -\frac 12(\sin 2\theta) \tau (P_{>0}(A)+P_{<0}(A))+ P_{\sigma_u}.$$
Furthermore, we fix $u_0$ and
define the unitary operator $V_u\in\cl(H)$ by
\begin{equation}
V_u|K^+:= U_{u}^* U_{u_0}, \quad V_u|K^-\oplus(\ker A)^\perp:=I.\mylabel{G7-2S.7}
\end{equation}
Then we choose a smooth family of self--adjoint operators, $T_u$, such
that
\begin{equation} V_u= e^{i T_u},\quad T_{u_0}=0,\quad T_u|K^-\oplus(\ker A)^\perp=0.
    \mylabel{G7-2S.8}
\end{equation}
It follows that
$$ V_u P_{u_0}(\theta) V_u^*=P_u(\theta)$$
and one checks that $(D_{\theta,\sigma_u})_u$
is a deformation of {\rm APS} type. Since $T_u'$ is an operator
of finite rank, we have
$$ \res(\gamma i T'_u)=\res(\gamma i T_u' (\sgn A) P_u(\theta))=0.$$
We deduce from Theorem \plref{S6-2.5}
$$\frac{d}{du} \Res_1 \eta(D_{\theta,\sigma_u};0)=0,$$
and
\begin{epeqnarray}
    \frac{d}{du} \ovl{\eta}(D_{\theta,\sigma_u})&=&
    \frac{1}{2\pi} a_{00}(A,\gamma i T_u')\\
      &=&\frac{i}{2\pi}\lim_{u\to 0} \tr_{H}[\gamma T_u'e^{-tA^2}]\\
      &=&\frac{i}{2\pi}\tr_{K^+}[\gamma T_u']\\
     &=&\frac{1}{2\pi i} \tr_{K^+}[U_u^{-1}\frac{d}{du} U_u].
\end{epeqnarray}

\vspace{1em}\noindent
Next we deal with the deformation $(D_{\theta,\sigma})_{|\theta|<\pi/2}$:

\begin{prop}{S7-2S.2} $\Res_1\eta(D_{\theta,\sigma};0)$ is independent
of $\theta$ and
\begin{eqnarray*}
   \frac{d}{d\theta}\ovl{\eta}(D_{\theta,\sigma})&=&
   \frac{1}{2\pi} a_{00}(A,\gamma (\sgn A) \tau)\\
      &=& \frac{1}{2\pi} \LIM_{t\to 0} \tr_{H}\big[
         \gamma (\sgn A) \tau e^{-t A^2}\big].
\end{eqnarray*}
Here $\DST \LIM_{t\to 0}$ is a
common notation for the constant term in the
asymptotic expansion as $t\to 0$.\revmark
\end{prop}
\proof In view of \myref{G7-2S.6} we put
$$T(\theta):=-i (\sgn A) \tau \theta.$$
Then one checks that \myref{G5-1.31a}--\myref{G5-2.28} and \myref{G7-2S.0} are satisfied.
We want to apply Corollary \plref{S4-1.6} to compute $\frac{d}{d\theta} \ovl{\eta}(D_{\theta,\sigma})$.
Since $\res$ vanishes on operators of finite rank we may replace
$$\gamma i T'(\theta)=\gamma (\sgn A) \tau$$
by
$$\gamma((\sgn A)+\sigma) \tau$$
in the assumptions of Corollary \plref{S4-1.6}. 
Since
$$(\gamma((\sgn A)+\sigma) \tau)^2=I$$
we infer from Lemma \plref{S4-1.4} that
$\res(\gamma((\sgn A)+\sigma) \tau)=0$. Thus 
$\Res_1 \eta(D_{\theta,\sigma};0)$
is independent of $\theta$ and 
\begin{epeqnarray}
    \frac{d}{d\theta} \ovl{\eta}(D_{\theta,\sigma})&=&
     \frac{1}{2\pi} a_{00}(A,\gamma i T'(\theta))\\
        &=& \frac{1}{2\pi} a_{00}(A,\gamma (\sgn A)\tau)\\
        &=&\frac{1}{2\pi}\LIMZ
           \tr_{H}
           \big[\gamma (\sgn A) \tau e^{-t A^2}\big].
\end{epeqnarray}

\vspace{1em}\noindent
Finally, we present the gluing law. In this situation
\myref{G1-1.6}--\myref{G5-2.17} we have yet another structure:
namely, introducing (with same notation as in \myref{G1-1.8}, \myref{G1-1.10})
$$\mu:=\mat{0}{1}{-1}{0}$$
we see that
$$\mu^2=-I,\quad \mu\tau+\tau\mu=\tau\tilde\gamma+\tilde\gamma\mu=
   \mu\tilde A+\tilde A\mu=0.$$
This observation leads to

\begin{theorem}{S7-2S.3} {\rm (Gluing Law)} Consider the deformation
of {\rm APS} type, $(D_{\theta,\sigma})_{|\theta|<\pi/2}$, introduced
in \myref{G7-2S.1}--\myref{G7-2S.6}. If there exists a unitary
classical pseudodifferential operator $\mu:\cinf{S_N}\to
\cinf{S_N}$ satisfying
\begin{equation}
   \mu^2=-I,\quad \mu\tau+\tau\mu=\tau\gamma+\gamma\mu=
   \mu A+A\mu=0
   \mylabel{G7-2S.9}
\end{equation}
then
$$\frac{d}{d\theta} \ovl{\eta}(D_{\theta,\sigma})=0.$$
\end{theorem}
\proof In view of \myref{G7-2S.9} we have
$$  \mu \gamma (\sgn A) \tau+\gamma (\sgn A) \tau \mu=0,$$
hence
   $$\tr_{H}\big[ \gamma (\sgn A) \tau e^{-t A^2}\big]=0.$$
In particular $a_{00}(A,\gamma(\sgn A) \tau)=0$ and, by
Proposition \plref{S7-2S.2}, we reach the conclusion.\endproof

Naming Theorem \plref{S7-2S.3} the ''gluing law''
calls for an
explanation: we briefly explain how the usual
gluing law for the $\eta$--invariant
follows from Theorem \plref{S7-2S.3}.

We consider again the situation \myref{G1-1.6}--\myref{G5-2.17}.
Then we have
\begin{eqnarray}
      \tilde K^\pm&:=&K^\pm(\tilde A)=
           \ker \tilde A\cap \ker (\tilde \gamma \mp i)\nonumber\\
     &=&K^\pm(A)\oplus K^\mp(A)\nonumber\\
   &=:& K^\pm\oplus K^\mp,
\end{eqnarray}
i. e. $\tilde K^\pm$ is canonically isomorphic to $K^+\oplus K^-$ and
we will use this identification in the sequel.
As in (\plref{G7-2S.3},\plref{G7-2S.4}) we write involutions
$\sigma$ of $\ker \tilde A$ with $\tilde \gamma\sigma+\sigma
\tilde\gamma=0$ in the
form
\begin{equation}
   \sigma(T)=\mat{0}{T^*}{T}{0},
\end{equation}
where $T:\tilde K^+\longrightarrow \tilde K^-$ is an isometry.
The isometry corresponding to the distinguished involution
\begin{equation}
  \sigma:=-\tau|\ker\tilde A=-\mat{0}{1}{1}{0}\otimes I_H
\end{equation}
in \myref{Grev-3.8d} therefore corresponds to the isometry
$-I:K^+\oplus K^-\longrightarrow K^+\oplus K^-.$

Now let $\ovl{\eta}(D,\mcut,T)$  be the $\ovl{\eta}$--invariant
of the operator $\tilde D$ with boundary condition given by
\begin{equation}
    \tilde P_T:=P_>(\tilde A)\oplus P_\sigma(T).
\end{equation}
Putting 
\begin{equation}
\tilde P_T(\theta):= 
\cos^2\theta P_{>0}(\tilde A)+\sin^2\theta P_{<0}(\tilde A)
    -\frac 12(\sin 2\theta) \tau (P_{>0}(\tilde A)+P_{<0}(\tilde A))+
    P_{\sigma(T)}
\end{equation}
as in \myref{G6-2.10b} we obtain a deformation of APS type
$(\tilde D_{\theta,\sigma(T)})_{|\theta|<\pi/2}$.
Then,
$\ovl{\eta}(\tilde D_{\theta,\sigma(T)})$ is independent of $\theta$
by Theorem \plref{S7-2S.3}.

For $T=-I$ and $\theta=\pi/4$ the boundary condition is the continuous
transmission boundary condition, hence $\ovl{\eta}(\tilde D_{\pi/4,\sigma(-I)})
\equiv\ovl{\eta}(D,M)\,\mod \Z$, the $\ovl{\eta}$--invariant of the closure
of $D$ on $L^2(S)$. Furthermore, for $\theta=0$ we obtain
\begin{equation}
\ovl{\eta}(\tilde D_{0,\sigma(T)})\equiv \ovl{\eta}(D,\mcut,T)\,\mod\Z.
\end{equation}
Thus, for $T=-I$ we have proved
\begin{equation}
\ovl{\eta}(D,\mcut,-I)\equiv\ovl{\eta}(D,M)\,\mod \Z.
\end{equation}

For an arbitrary isometry $T:K^+\oplus K^-\longrightarrow
K^+\oplus K^-$ we choose a smooth path of isometries
$(T_u)_{0\le u\le 1}$ with $T_0=-I, T_1=T$ and apply Proposition
\plref{S7-2S.1} to $\tilde D_{0,\sigma(T_u)}$. Then
\begin{equation}
\frac{d}{du} \ovl{\eta}(D,\mcut,T_u)=\frac{1}{2\pi i} \tr_{K^+\oplus K^-}
(T_u^{-1} \frac{d}{du} T_u)=\frac{1}{2\pi i}\frac{d}{du} \log\det(T_u),
\end{equation}
and hence
\begin{equation}
\ovl{\eta}(D,\mcut,T)\equiv\ovl{\eta}(D,\mcut,-I)+\frac{1}{2\pi i}
  \log \det(T)-\frac 12(\dim K^++\dim K^-)\,\mod\, \Z.
\end{equation}

This can be written more nicely in terms of the $\tau$--invariant
\myref{G4-1.21}. Namely,
\begin{equation}
    \tau(D,\mcut,T)=(-1)^{\dim K^++\dim K^-} \det(T) \tau(D,M).
\end{equation}
Note that $A_+:=A|\ker(\gamma-i)$ is a Fredholm operator between
$\ker(\gamma-i)$ and $\ker(\gamma+i)$ with $\ind A_+=\dim K^+-\dim K^-$
and hence we end up with the gluing law in the version of 
Dai and Freed \cite[Prop. 4.5]{DF}
\begin{theorem}{Srev-glu5}
 $\displaystyle   \tau(D,\mcut,T)=(-1)^{\ind A_+} \det(T) \tau(D,M).$
\end{theorem}
Actually, this result is slightly more general than loc. cit. since
Dai and Freed deal with Dirac operators on spin manifolds.

In the special case that the hypersurface $N$ separates $M$ into
two components $M_\pm$ such that $M=M_-\cup_N M_+$, the index of 
$A_+$ vanishes by the cobordism theorem (cf. the discussion on
page \pageref{cobordism}). Hence we can choose isometries
$T_+:K^+\longrightarrow K^-$, $T_-:K^-\longrightarrow K^+$ and put
\begin{equation}
T:=\mat{0}{T_+}{T_-}{0}.
\end{equation}

Then 
\begin{equation}
    \tilde D_{0,\sigma(T)}=D^+_{0,\sigma(T_+)}\oplus D^-_{0,\sigma(T_-)},
\end{equation}
where  $D^\pm_{0,\sigma(T_\pm)}$ is the operator $D$ on the manifold
$M_\pm$ with boundary condition given by $P_>(A)\oplus P_{\sigma(T_+)}$
resp. $P_>(-A)\oplus P_{\sigma(T_-)}$. Denoting their respective
$\ovl{\eta}$--invariants by $\ovl{\eta}(D,M_\pm,T_\pm)$ we obtain
the gluing formula for the $\eta$--invariant
\begin{equation}
\ovl{\eta}(D,M)\equiv \ovl{\eta}(D,M_+,T_+)+
      \ovl{\eta}(D,M_-,T_-)+\frac{1}{2\pi i} \log \det(-T_1 T_2)\,\mod\Z,
\end{equation}
or, in multiplicative notation,
\begin{equation}
\tau(D,M)=\det(-T_1 T_2) \tau(D,M_+)\tau(D,M_-).
\end{equation}

As explained in \cite[Sec.1]{Bunke} $\det(-T_1 T_2)$ is related to the
Maslov index of the corresponding Lagrangian subspaces defined by
$L_j:=\ker\sigma(T_j), j=1,2$. Namely, putting
\cite[Theorem 2.1]{LeschW}, \cite[Def. 1.3]{Bunke}
$$m(L_1,L_2)=-\frac 1\pi \sum_{\begin{array}{c}\SST \beta\in (-\pi,\pi)\\
  \SST     e^{i\beta}\in\spec(-T_1T_2)\end{array}} \beta,$$
then
$$m(L_1,L_2)=\int_K \tau(kL,L_1,L_2)dk.$$
Here, $K$ is the stabilizer of $\tilde \gamma$ in the symplectic group,
$L$ is an arbitrary Lagrangian subspace,
and $\tau$ is the Maslov triple index (cf. \cite[Sec. 2]{Bunke} for details).

Summing up we can state the gluing law as follows:

\begin{theorem}{Srev-glu6}
$$\ovl{\eta}(D,M)\equiv \ovl{\eta}(D,M_+,T_+)+
      \ovl{\eta}(D,M_-,T_-)+\frac 12 
  m(\ker\sigma(T_1),\ker\sigma(T_2))\,\mod\Z.$$
\end{theorem}

\vspace*{1em}
Our last comment concerns the residue at $0$ of the $\eta$--function.
We expect that in general the residue in \myref{G6-2.36}
will not vanish.
In the cutting case, however, there is no pole:

\begin{theorem}{S5-2.4} If $(\dth)_{|\theta|<\pi/2}$ arises from cutting
$M$ along a compact hypersurface {\rm (}as
explained in \myref{G1-1.6}--\myref{G5-2.17}{\rm )}
then $\eta(\dth;s)$ is regular
at $s=0$, for all $\theta$.
\end{theorem}

\proof  By Proposition \plref{S7-2S.2}, $\Res_1 \eta(D_\theta;0)$ is
independent of $\theta$, hence
$$\Res_1 \eta(D_\theta;0)=\Res_1\eta(D_{\pi/4};0)=0$$
since the $\eta$--function of a self--adjoint elliptic differential
operator on a compact manifold is regular at $0$
\cite[Sec. \mls 3.8]{Gi}.\endproof

\def\ti{{\widetilde{\rm I}}}
\def\tii{{\widetilde{\rm II}}}
\def\tiii{{\widetilde{\rm III}}}
\def\nti{{\rm I}}
\def\ntii{{\rm II}}
\def\ntiii{{\rm III}}

\section{Proofs}\label{Sec3}

We now prove the statements used in the previous section.

\par\medbreak\noindent
{\bf Proof of Proposition \ref{S5-2.1}}\quad
We consider $\tilde\dth$ first. Let $u\in\cd(\tilde\dth^*)$ satisfy
$$\tilde D^*_\theta u=\pm \sqrt{-1} u.$$
This implies, for $v\in\tilde\cd_\theta$, that
$$(\tilde D_\theta v,u)=\mp \sqrt{-1} (v,u).$$
Then a standard regularity argument shows that 
$u\in C(\R_+,L^2(S_N))$ with
$$P(\theta)u(0)=0,$$
by \myref{G5-1.31a}. Choosing $\phi\in\cinfz{\R}$ with
$\phi=1$ near $0$ we put $\phi_N(x):=\phi(x/N)$ and obtain
$\phi_N^2 u\in\tilde\cd_\theta$. Consequently, we
find that
$$\pm\sqrt{-1} \|u\|^2=\lim_{N\to\infty} (\tilde\dth\phi_N^2 u,u)
    =\lim_{N\to\infty} (u,\tilde\dth\phi_N^2 u)\in\R,$$
hence $u=0$.

For $\dth$, we appeal to the localization principle for
deficiency indices derived in 
\cite[Thm. \mls 2.1]{Lesch1} (cf. also  \cite[Chapter IV]{Lesch2}).\endproof

In what follows it will be crucial that we can give an
explicit formula for the operator heat kernel of $\tilde\dth$.
It is the operator analogue of a formula derived by Sommerfeld
\cite[p. \mls 61]{Sommer}.

\begin{theorem}{S5-3.1}  We have for $t,x,y>0$
\begin{eqnarray}
  e^{-t\tilde\dth^2}(x,y)&=&(4\pi t)^{-1/2} 
      \left(e^{-(x-y)^2/4t}+(I-2P(\theta))e^{-(x+y)^2/4t}\right) e^{-tA^2}
         \nonumber\\
    &&\quad + (\pi t)^{-1/2}(I-P(\theta))
      \int_0^\infty e^{-(x+y+z)^2/4t} \tilde A(\theta) e^{\tilde A(\theta)
       z-t A^2} dz,\mylabel{G5-3.1}
\end{eqnarray}    
where $\tilde A(\theta):=(I-P(\theta))A(I-P(\theta))$.
\end{theorem}
\proof The point is the convergence of the integral in
\myref{G5-3.1}. Note that $P(\theta)$ commutes with $|A|$ by
\myref{G5-1.31b} and the discreteness of $A$. Thus from
\myref{G5-1.31a}, \myref{G1-1.3}, and \myref{G5-2.25}
\begin{eqnarray*}
   \tilde A(\theta)&=& \gamma P(\theta) \gamma^* A \gamma P(\theta)\gamma^*
                        =-\gamma P(\theta) A P(\theta)\gamma^*\\
                   &=&-a(\theta) \gamma |A| P(\theta)\gamma^*\\
                   &=& -a(\theta) |A| (I-P(\theta)).
\end{eqnarray*}
In particular, $\tilde A(\theta)$ commutes with $(I-P(\theta))$ so
$$ \tilde A(\theta) e^{\tilde A(\theta)z-t A^2}
   = -a(\theta) |A| (I-P(\theta))e^{-a(\theta)|A|z-t A^2}.$$
Introducing $a_-(\theta):= -\min\{0,a(\theta)\}\in [0,1)$ we find
$$-a(\theta) |A| z\le a_-(\theta) \left(\frac{z^2}{4t}+A^2t\right),$$
and
\begin{equation}
   0\le |A|(I-P(\theta)) e^{\tilde A(\theta)z- tA^2}\le
       |A| (I-P(\theta))e^{a_-(\theta) z^2/4t} e^{-(1-a_-(\theta))t A^2}.
       \mylabel{G5-3.2}
\end{equation}
This implies that the integral converges in the trace norm of
$L^2(S_N)$.

Now pick $u\in\cinfz{(0,\infty),L^2(S_N)}$ and form
$$Q_t u(x):=\int_0^\infty Q_t(x,y) u(y) dy$$
where $Q_t$ denotes the right hand side of \myref{G5-3.1}.
Then it is a routine matter to check that we have
\begin{equation}\begin{array}{l}
  \DST Q_tu\in C^1((0,\infty),\cd(\tilde D^*))\cap
      C(\R_+,\ch),\\
  \DST (\pl_t + (\tilde D^*)^2 Q_t u(x)=0,\quad t,x>0,\\
  \DST \lim_{t\to 0+} Q_t u(x)=u(x).
  \end{array}
\end{equation}
Hence it remains to verify the boundary conditions. Clearly,
$$P(\theta)Q_t(x,y)=(4\pi t)^{-1/2}
    \left(e^{-(x-y)^2/4t}-e^{-(x+y)^2/4t}\right)P(\theta) e^{-tA^2}
    \marrow_{x\to 0+} 0,$$
and the same holds for $P(\theta) Q_t u(x)$ and $A P(\theta) Q_tu(x)$,
by dominated convergence.
This implies
$$Q_t u\in \cd(\tilde D_\theta).$$
We finally have to show that
\begin{eqnarray*}
    0&=& \lim_{x\to 0+} P(\theta) \gamma (\pl_x+A)Q_t u(x)\\
     &=&\lim_{x\to 0+}\gamma (I-P(\theta)) (\pl_x+A)Q_tu(x)\\
     &=&\lim_{x\to 0+}\left\{\gamma(\pl_x+\tilde A(\theta))
         (I-P(\theta))Q_tu(x)+\gamma (I-P(\theta))A P(\theta)Q_tu(x)
         \right\}\\
     &=&\lim_{x\to 0+}\gamma(\pl_x+\tilde A(\theta))
         (I-P(\theta))Q_tu(x).
\end{eqnarray*}
An easy calculation shows that
\begin{eqnarray*}
   &&(\pl_x +\tilde A(\theta))(I-P(\theta)) Q_t(x,y)\\
   &=& (4\pi t)^{-1/2}\left\{ e^{-(x-y)^2/4t} {\TST(\frac{y-x}{2t} + \tilde A(\theta))}
       +  e^{-(x+y)^2/4t} {\TST(-\frac{y+x}{2t} +\tilde A(\theta))}\right\}
        (I-P(\theta)) e^{-t A^2}\\
   &&\quad - (\pi t)^{-1/2} e^{-(x+y)^2/4t}\tilde A(\theta) (I-P(\theta))
      e^{-tA^2} \marrow_{x\to 0+} 0.
\end{eqnarray*}
Then the proof is completed as above.\endproof

\par\medbreak\noindent
{\bf Proof of Proposition \ref{S5-2.2}}\quad
We propose to show that, for $u\in \cd(D_\theta^k)$ with
$k>m/2$, we have the estimate
\begin{equation}
     |u(p)|\le C (1-a_-(\theta))^{-k-1/2} 
     \big(\|u\|_{L^2(S)}+\|D_\theta^k u\|_{L^2(S)}\big).
     \mylabel{G5-3.4}
\end{equation}
As explained in \cite{Lesch3} (cf. also \cite[Sec. 1.4]{Lesch2}),
this estimate implies the Hilbert--Schmidt property of suitable
functions of $D_\theta$ and, in particular, the assertion of the
proposition.

To prove \myref{G5-3.4}, it is clearly enough to assume that
$\supp u\subset U$, and we are reduced to proving the analogue
of \myref{G5-3.4} for $\tilde\dth$ if $\supp u\subset [0,1)$.
To do so, we write for
$q\in N$
\begin{eqnarray}
   u(x)(q)&=&(\tilde \dth^2+1)^{-j}(\tilde \dth^2+1)^j u(x)(q)\nonumber\\
   &=& \frac{1}{\Gamma(j)}
       \int_0^\infty e^{-t} t^{j-1}\int_0^\infty
        e^{-t \tilde \dth^2}(x,y)(\tilde \dth^2+1)^j u(y) dy dt (q).
\end{eqnarray}
From the ellipticity of $A$ we get for $k> (m-1)/2$
$$|u(x) (q)| \le C_k \| (A^2+1)^k u(x)\|_{L^2(S_N)},$$
hence, with $j=k+1/2+\eps$, $\eps>0$,
\begin{eqnarray}
  |u(x)(q)|&\le& C_k \int_0^\infty e^{-t} t^{k-1/2+\eps}\int_0^\infty
     \|(A^2+1)^k e^{-t\tilde\dth^2}(x,y)\cdot\nonumber\\
   &&\quad\quad\cdot((\tilde\dth^2+1)^{k+1/2+\eps}u(y))\|_{L^2(S_N)}dydt
     \mylabel{G5-3.6}
\end{eqnarray}

From \myref{G5-3.1} and \myref{G5-3.2} we derive the norm estimate
\begin{eqnarray}
   &&\|(A^2+1)^k e^{-t\tilde\dth^2}(x,y)\|_{\cl(L^2(S_N))}\nonumber\\
   &\le&C_k (1-a_-(\theta))^{-k-1}\;t^{-k-1/2} \left(
        e^{-(x-y)^2/4t}+ e^{-(x+y)^2/4t}\right).
        \mylabel{G5-3.7}
\end{eqnarray}
Using \myref{G5-3.7} and the Cauchy--Schwarz inequality in \myref{G5-3.6}
we obtain the result.\endproof


\bigskip\noindent
{\bf Proof of Lemma \ref{S6-3.2}}\quad
An integration by parts gives
\begin{eqnarray}
   F_a(x)&=&-\frac{1}{2a} \int_0^\infty \erfc(z) \frac{\pl}{\pl z}
           \big(e^{-2a xz-x^2}\big)dz\nonumber\\
      &=& \frac{1}{2a} e^{-x^2} -\frac{1}{a\sqrt{\pi}}\int_0^\infty
           e^{-(2a xz +x^2+z^2)} dz\nonumber\\
      &=:& \frac{1}{2a}\Big( G(x)-\tilde F_a(x)\Big).
      \mylabel{G6-L32.6}
\end{eqnarray}
Clearly,
\begin{equation}
   {\cm} G(w)= \frac 12 \Gamma(w/2).
      \mylabel{G6-L32.7}
\end{equation}
To determine ${\cm} \tilde F_a$, we observe that
$$\tilde F_a(x) = e^{-(1-a^2) x^2} \erfc(ax)$$
and derive a differential equation in $a$. In fact, for $\Re w>0$,
$0<|a|<1$,
\begin{eqnarray}
    \frac{\pl}{\pl a}(1-a^2)^{w/2} {\cm} \tilde F_a(w) &=&
         \frac{\pl}{\pl a} \int_0^\infty x^{w-1} e^{-x^2}\erfc(
            \frac{a}{\sqrt{1-a^2}}x) dx\nonumber\\
        &=&-\frac{2}{\sqrt{\pi}} \int_0^\infty x^w
             e^{-x^2/(1-a^2)} dx \;(1-a^2)^{-3/2}\nonumber\\
        &=& -\frac{1}{\sqrt{\pi}} \Gamma((w+1)/2)(1-a^2)^{w/2-1}
           \mylabel{G6-L32.8}
\end{eqnarray}
The initial condition at $a=0$ is
\begin{equation}
     {\cm} \tilde F_0(w)= \frac 12 \Gamma(w/2).
\end{equation}
The solution of this initial value problem is, for $|a|<1$,
$${\cm} \tilde F_a(w)
= (1-a^2)^{-w/2} \Big( \frac 12 \Gamma(w/2)-\frac{1}{\sqrt{\pi}}
   \Gamma((w+1)/2) \int_0^a (1-t^2)^{w/2-1}dt \Big),$$
hence
\begin{eqnarray}
  {\cm} F_a(w)&=&\frac{1}{2a}\Big[(1-(1-a^2)^{-w/2})\frac 12\Gamma(w/2)\nonumber\\
      &&+\frac{1}{\sqrt{\pi}}(1-a^2)^{-w/2}
         \int_0^a (1-t^2)^{w/2-1} dt \;\Gamma((w+1)/2)\Big].
          \mylabel{G6-L32.9}
\end{eqnarray}
Furthermore,
\begin{eqnarray}
   {\cm} \tilde F_1(w)&=& \int_0^\infty x^{w-1} \erfc(x) dx\nonumber\\
      &=&\frac{2}{w\sqrt{\pi}} \int_0^\infty x^w e^{-x^2} dx
         =\frac{1}{w\sqrt{\pi}} \Gamma((w+1)/2),
\end{eqnarray}
thus
\begin{equation}
{\cm} F_1(w)=\frac 14\Big[\Gamma(w/2)-\frac{2}{w\sqrt{\pi}}
            \Gamma((w+1)/2)\Big].
\end{equation}
The poles and residues of ${\cm} F_a$ can now easily be calculated
in terms of the poles and residues of the $\Gamma$--function.
\endproof

\noindent
We turn to the

\medskip

\noindent
{\bf Proof of Theorem \ref{S1-2.7}}
We choose $\phi\in\cinfz{-1,1}$ with $\phi=1$
near $0$. Then, from \cite[Lemma 1.9.1]{Gi}
(cf. Remark 2) after Lemma \plref{S6-1.2.5})
we obtain the asymptotic expansion, for $l=0,1$,
\begin{equation}
   \tr_{L^2(S)}[(1-\psi_\phi)D_\theta^l e^{-tD_\theta^2}]\sim_{t\to 0+}
     \sum_{j=0}^\infty a_j(\phi;\theta,l)\;t^{j-m/2}.
     \mylabel{G2-2.29}
\end{equation}
The coefficients can be computed locally in terms of the
natural extension of $D$ to the metric double,
$\tilde M$, of $M$, and $\psi_\phi$.

Thus, since $e^{-t\tilde D^2}$ can serve as a parametrix for
$D_\theta^2$ we obtain from 
\cite[Theorem 2.10 and Prop. 3.4]{Lesch3} (note that
what is called there the ''singular elliptic estimate''
was proved in \myref{G5-3.4})
that
\begin{equation}
\tr_{L^2(S)}[\psi_\phi D_\theta^l e^{-tD_\theta^2}]
   \sim_{t\to 0+} \tr_\ch[\phi\tilde D_\theta^l
    e^{-t\tilde D_\theta^2}],
    \mylabel{G6-3.9}
\end{equation}    
and it is enough to expand the right hand side of \myref{G6-3.9}
for $l=0,1$.

Consider $l=0$ first. We obtain from the explicit formula
\myref{G5-3.1} and the Trace Lemma \cite[Appendix]{BS}
that
\begin{eqnarray}
   \tr_\ch[\phi
    e^{-t\tilde D_\theta^2}]&=&
      \int_0^\infty \phi(x) (4\pi t)^{-1/2} \tr_H [e^{-tA^2}] dx\nonumber\\
      && + \int_0^\infty \phi(x) e^{-x^2/t} (4\pi t)^{-1/2}
         \tr_H[ (I-2P(\theta)) e^{-tA^2}] dx\nonumber\\
      &&-a(\theta) \int_0^\infty \int_0^\infty
        \phi(x) e^{-(2x+z)^2/4t} (\pi t)^{-1/2} 
          \tr_H[P(\theta) |A| e^{-a(\theta) |A| z -t A^2}] dz dx\nonumber\\
   &=:&\nti(t) +\ntii(t)+\ntiii(t).\mylabel{G6-3.10}
\end{eqnarray}    

Since $A$ is elliptic on $S_N$ we have for the first term
\begin{equation}
  \nti(t)\sim_{t\to 0+} (4\pi t)^{-1/2} \int_0^\infty \phi(x) dx
     \sum_{j=0}^\infty b_j(A^2) t^{j-(m-1)/2}.
     \mylabel{G6-3.11}
\end{equation}
Next, as an easy consequence of \myref{G5-1.31a} we see that
\begin{equation}
  \ntii(t)\equiv 0.\mylabel{G6-3.12}
\end{equation}

For $\ntiii(t)$, we write, with
$$c(\gl):=\dim\ker(|A|-\gl)=2\, \tr_{\ker(|A|-\gl)}(P(\theta)),$$
\begin{eqnarray}
   \ntiii(t)&=& -a(\theta) \int_0^\infty\int_0^\infty \phi(x\sqrt{t})
        \frac{1}{\sqrt{\pi}} e^{-(x+z)^2}
        \sum_{\gl\in\spec |A|\setminus\{0\}}
        c(\gl) \sqrt{t} \gl e^{-2a(\theta) \sqrt{t}\gl z-t\gl^2}dzdx
        \nonumber\\
     &\sim_{t\to 0+    }& -\frac{a(\theta)}{2}
       \int_0^\infty \erfc(z) \sum_{\gl\in\spec |A|\setminus\{0\}}
        c(\gl) \sqrt{t} \gl e^{-2a(\theta) \sqrt{t}\gl z-t\gl^2}dz\nonumber\\
    &=& -\frac{a(\theta)}{2} \sum_{\gl\in\spec |A|\setminus\{0\}}
        c(\gl) F_{a(\theta)}(\sqrt{t}\gl)\nonumber\\
    &=& -\frac{a(\theta)}{2} \sum_{\gl\in\spec |A|\setminus\{0\}}
       c(\gl) \frac{1}{2\pi i} \int_{\Re w=c>>0}
       t^{-w/2} \gl^{-w} {\cm} F_{a(\theta)}(w) dw\nonumber\\
    &=&-\frac{a(\theta)}{4\pi i}\int_{\Re w=c}
        t^{-w/2} \zeta_{A^2}(w/2) {\cm} F_{a(\theta)}(w) dw.
        \mylabel{G6-3.13}
\end{eqnarray}
We now collect the various contributions. First, replacing
$\phi$ by $\phi_\eps$, $\phi_\eps(x):=\phi(x/\eps)$,
and letting $\eps\to 0$ we obtain from \myref{G2-2.29} and
\myref{G6-3.11} a contribution
\alpheqn
\begin{equation}
  {\rm \tilde I}(t)\sim_{t\to 0+} \sum_{j=0}^\infty
         a_j(\theta,0) \;t^{j-m/2},
         \mylabel{G6-3.14a}
\end{equation}
where
$$a_j(\theta,0)= \int_M \tilde u_j(\theta,0),$$
with $\tilde u_j$ a local density computed for the natural
extension of $D$ to the double, $\tilde M$, of $M$.

The remaining contribution, $\ntiii(t)$, can be evaluated by the
Residue Theorem since the integrand decays in
vertical strips with bounded real part (by Lemma \plref{S6-3.2},
Lemma \plref{S6-1.2.5}, and \myref{G7-1.19}). Thus we find
(using e.g. the description of the singularities of
$\zeta_{A^2}$ in \cite[Lemma 2.1]{BL})
\begin{eqnarray}
\ntiii(t)&= &-\frac {a(\theta)}{2}
   \sum_{w\in\C} \Res_1 (t^{-w/2}\zeta_{A^2}(w/2) {\cm} F_{a(\theta)}(w))
   \nonumber\\
   &\sim_{t\to 0+}& \frac{a(\theta)}{2}
      \sum_{j=0}^\infty t^{j-n/2}\Big\{\log t\, \Res_1 \zeta_{A^2}(n/2-j)
        \Res_1 {\cm} F_{a(\theta)}(n-2j)\nonumber\\
  &&-2\Res_1 \zeta_{A^2}(n/2-j) \Res_0 {\cm} F_{a(\theta)} (n-2j)\Big\}\nonumber\\
  &&-\frac{a(\theta)}{2} \sum_{j=0}^\infty
      t^{j/2} \Res_0 \zeta_{A^2}(-j/2) \Res_1 {\cm} F_{a(\theta)} (-j).
      \mylabel{G6-3.14b}
\end{eqnarray}
\reseteqn

From this, we can read off our assertions on the structure
of the coefficients. First of all, the leading contribution
comes from \myref{G6-3.14a} only, as
$a_0 t^{-m/2}$, and so is computed as in the compact
case. Next, we observe that $\zeta_{A^2}$ has no poles
at the points $n/2-j$ for $j\ge n/2$ if $n$ is even. If
$n$ is odd, however, the $\log$--terms occur as can be seen
from Lemma \plref{S6-3.2}.
The coefficients of the terms in the first sum in
\myref{G6-3.14b} are computed from local densities
on $N$, whereas those in the second sum are, in general,
nonlocal.

Next we consider the case $l=1$. In view of \myref{G2-2.29}
and the previous analysis
it is enough to expand
\begin{equation}
   \int_0^\infty \phi(x) \tr_H[\gamma (\pl_x+A)e^{-t\tilde D_\theta^2}
   (x,x)]dx
   =:\ti(t)+\tii(t)+\tiii(t),
       \mylabel{G6-3.15}
\end{equation}
numbering again the contributions according to the three terms
in \myref{G5-3.1}. In view of \myref{G5-1.31a},
\myref{G5-1.31b}, \myref{G1-1.3}, and \myref{G7-2S.2} we find
\begin{equation}\begin{array}{l}
   \DST \tr_H[\gamma e^{-tA^2}]=\tr_H[\gamma P(\theta) e^{-tA^2}]=0\\[0.5em]
   \DST \tr_H[\gamma P(\theta)|A| e^{-a(\theta)|A|z-tA^2}]=0,
   \end{array}
   \mylabel{G6-3.17}
\end{equation}
and thus
\begin{equation}
   \tr_H[\gamma \pl_x e^{-t \tilde D_\theta^2}(x,x)]=0.
   \mylabel{G6-3.18}
\end{equation}
Again from \myref{G1-1.3} we conclude
\begin{equation}
   \tr_H[\gamma A e^{-tA^2}]=0,
   \mylabel{G6-3.19}
\end{equation}
which implies
\begin{equation}
     \ti(t)\equiv 0.
        \mylabel{G6-3.20}
\end{equation}
Furthermore,
\begin{eqnarray}
    \tii(t)&=&  (4\pi t)^{-1/2}\int_0^\infty \phi(x) e^{-x^2/t}
         \tr_H[\gamma A (I-2P(\theta)) e^{-tA^2}] dx\nonumber\\
         &=&(4\pi)^{-1/2}\int_0^\infty \phi(x\sqrt{t}) e^{-x^2}dx\,
         \tr_H[\gamma A (I-2P(\theta)) e^{-tA^2}]\nonumber\\
         &\sim_{t\to 0+}& \frac 14 \tr_H[\gamma A (I-2P(\theta)) e^{-tA^2}]\nonumber\\
         &=& -\frac 12\tr_H[\gamma A P(\theta) e^{-tA^2}].
         \mylabel{G6-3.21}
\end{eqnarray}
Finally, we note that, using again \myref{G5-1.31a}
and \myref{G1-1.3},
\begin{eqnarray}
     \tr_H[\gamma A (I-P(\theta)) \tilde A(\theta) e^{\tilde A(\theta)z-tA^2}]
     &=& a(\theta)\tr_H[\gamma A P(\theta)|A| e^{-a(\theta)|A|z-tA^2}],
     \mylabel{G6-3.22}
\end{eqnarray}
and so, as in \myref{G6-3.13}, with
$d(\gl):= \tr_{\ker(|A|-\gl)}[\gamma A P(\theta)]$,
\begin{eqnarray}
   \tiii(t)&=&a(\theta) \int_0^\infty\int_0^\infty \phi(x)
         e^{-(2x+z)^2/4t}(\pi t)^{-1/2}
          \tr_H[\gamma A P(\theta)|A| e^{-a(\theta)|A|z-tA^2}]
             dzdx\nonumber\\
     &=& a(\theta) \int_0^\infty\int_0^\infty \phi(x\sqrt{t})
        \frac{2}{\sqrt{\pi}} e^{-(x+z)^2}
        \sum_{\gl\in\spec |A|\setminus\{0\}}
        d(\gl) \sqrt{t} \gl e^{-2a(\theta) \sqrt{t}\gl z-t\gl^2}dzdx
        \nonumber\\
    &\sim_{t\to 0+}& a(\theta) \sum_{\gl\in\spec |A|\setminus\{0\}}
        d(\gl) F_{a(\theta)}(\sqrt{t}\gl)\nonumber\\
    &=&\frac{a(\theta)}{2\pi i}\int_{\Re w=c}
        t^{-w/2} \eta(A,\gamma A P(\theta);w-1) {\cm} F_{a(\theta)}(w) dw.
        \mylabel{G6-3.23}
\end{eqnarray}
Combining our computations, we see that the terms local
on $\tilde M$ protrude from \myref{G2-2.29} as before.

We obtain the second contribution from \myref{G6-3.21}.
However, since $P(\theta)$ is a pseudodifferential
operator we now have to employ the general expansion
theorem for pseudodifferential operators \myref{G7-1.19}
\cite[Theorem 2.7]{GrubbSeeley}.
Namely
\alpheqn
\begin{eqnarray}
  \tii(t)&\sim_{t\to 0+}&
      -\frac 12 \tr_H[\gamma A P(\theta) e^{-tA^2}]\nonumber\\
       &\sim_{t\to 0+}&\sum_{j=0}^\infty
          c_j^1(\theta,1) t^{(j-m)/2}
          +\sum_{j=0}^\infty \Big(b_j^1(\theta,1) t^j \log t+
            d_j^1(\theta,1) t^j\Big).
            \mylabel{G6-3.24}
\end{eqnarray}
Here, $b_j^1,c_j^1$ are integrals of local densities over $N$ whereas the
$d_j^1(\theta,1)$ are, in general, nonlocal spectral invariants on $N$.

For the third contribution, we use again the estimate
\myref{G7-1.13} with $B=\gamma A P(\theta)$
(stemming from the fact that $P(\theta)$ is a pseudodifferential
operator) to obtain
 $$  \tiii(t)\sim_{t \to 0+}
   a(\theta)\sum_{w\in\C} \Res_1\Big(t^{-w/2} \eta(A,\gamma A P(\theta);w-1)
       {\cm} F_{a(\theta)}(w)\Big).$$
From the expansion \myref{G6-3.24} and Lemma \plref{S4-1.1}
one derives that $\eta(A,\gamma A P(\theta);w)$
is meromorphic in $\C$ with simple poles at the points $n-k, k\in \Z_+$.
Furthermore, the residues of the poles are integrals of local
densities over $N$. Thus
\begin{eqnarray}
   \tiii(t)&\sim_{t \to 0+}&
       -\frac{a(\theta)}{2}
       \sum_{j=0}^\infty t^{j/2} \log t\;
       \Res_1(\eta(A,\gamma A P(\theta);-j-1)\Res_1 {\cm} F_{a(\theta)}(-j)\nonumber\\
       &&+a(\theta) \sum_{j=0}^\infty t^{(j-m)/2}
          \Res_1(\eta(A,\gamma A P(\theta);m-j-1)\Res_0 {\cm} F_{a(\theta)}(m-j)
          \nonumber\\
       &&+a(\theta) \sum_{j=0}^\infty t^{j/2}
          \Res_0(\eta(A,\gamma A P(\theta);-j-1)\Res_1 {\cm} F_{a(\theta)}(-j).
          \mylabel{G6-3.26}
\end{eqnarray}
\reseteqn
The coefficients in the first and second sum are again local, like
$c_j^1$ in \myref{G6-3.24}, whereas those in the second sum are not.

It remains to compute the contribution to $t^{-1/2}$
from (\ref{G6-3.24},b). Using Lemma \plref{S4-1.1},
it turns out to be equal to
\begin{eqnarray*}
    &&-\frac 12 a_{-1/2,0}(A,\gamma A P(\theta))+
        a(\theta) \Res_0 {\cm} F_{a(\theta)}(1)\Res_1 \eta(A,\gamma A P(\theta);0)\\
    &=&\big(-\frac{\sqrt{\pi}}{4}+a(\theta) {\cm} F_{a(\theta)}(1))
        \Res_1 \eta(A,\gamma A P(\theta);0)\\
    &=&\big(-\frac{\sqrt{\pi}}{4}+a(\theta) {\cm} F_{a(\theta)}(1))
        \Res_1 \eta(A,\gamma (\sgn A) P(\theta);-1)\\
    &=&\big(-\frac{\sqrt{\pi}}{4}+a(\theta) {\cm} F_{a(\theta)}(1))
        \res (\gamma (\sgn A) P(\theta)),
\end{eqnarray*}
using \myref{G5-1.23} in the last step.\endproof

\bigskip\noindent
{\bf Proof of Theorem \ref{S6-2.5}}
We choose $\tilde \phi\in\cinfz{-1,1}$ with $\tilde \phi=1$ in a
neighborhood of $\supp \phi$, with
$\phi$ from \myref{G5-2.31}. Then, for $u\in \cd_0$, one
easily computes (writing $\tilde \psi=\psi_{\tilde \phi}$)
\begin{eqnarray}
    \check D_\theta \tilde\psi u&=& \Phi_\theta^*\Phi^* \gamma
       (\pl_x+A)\Phi \Phi_\theta \tilde \psi u\\
       &=:& \Phi^* i \gamma \phi' T(\theta)\Phi \tilde\psi u
         +\Phi^*\Phi_\theta^*\gamma A \Phi_\theta \Phi \tilde \psi u
         +v,
\end{eqnarray}
with $v$ independent of $\theta$, hence
$$\frac{d}{d\theta} \check D_\theta \tilde \psi u
  = \Phi_\theta^*\Phi^* i\gamma(\phi'T'(\theta)-2\phi T'(\theta)A)
      \Phi \Phi_\theta \tilde \psi u$$
and
\begin{equation}
\tr_{L^2(S)}\big[\frac{d}{d\theta} \check D_\theta e^{-t\check D_\theta^2}\big]
  = i\tr_{L^2(S)}\big[\gamma(\phi'T'(\theta)-2\phi T'(\theta)A)
      e^{-t D_\theta^2}\tilde\psi\big].
   \mylabel{G6-3.29}
\end{equation}
We can argue as in the proof of Theorem \ref{S1-2.7} to replace
$e^{-t D_\theta^2}$ by $e^{-t \tilde D_\theta^2}$, i.e.
\begin{eqnarray}
  &&   i\tr_{L^2(S)}\big[\gamma(\phi'T'(\theta)-2\phi T'(\theta)A)
      e^{-t D_\theta^2}\tilde\psi\big]\nonumber\\
  &\sim_{t\to 0+}&   i\tr_{L^2(S)}\big[\gamma(\phi'T'(\theta)-2\phi T'(\theta)A)
      e^{-t \tilde D_\theta^2}\tilde\psi\big].
   \mylabel{G6-3.30}
\end{eqnarray}
Again as in the proof of Theorem \ref{S1-2.7}, we obtain
twice three terms from plugging the kernel \myref{G5-3.1}
in \myref{G6-3.30}.

We start with
\begin{eqnarray}
     && i\tr_{L^2(S)}[\gamma\phi'T'(\theta)e^{-t\tilde D_\theta^2}]\nonumber\\
     &=&i\int_0^\infty \phi'(x) \tr_H[\gamma T'(\theta) e^{-t\tilde D_\theta^2}(x,x)]dx
        \nonumber\\
     &=:&\nti(t)+\ntii(t)+\ntiii(t).\mylabel{G6-3.31}
\end{eqnarray}
We find
\begin{eqnarray}
   \nti(t)&=& i (4\pi t)^{-1/2}  \int_0^\infty \phi'(x) dx \;\tr_H[\gamma T'(\theta)
           e^{-tA^2}]\nonumber\\
       &=&  -i (4\pi t)^{-1/2} \tr_H[\gamma T'(\theta)e^{-tA^2}].
       \mylabel{G6-3.32}
\end{eqnarray}
Since $\phi'$ is supported away from zero, it is easy to see
that
\begin{equation}
   \ntii(t)\sim_{t\to 0+} \ntiii(t)\sim_{t\to 0+} 0.
   \mylabel{G6-3.33}
\end{equation}

The second contribution is
\begin{eqnarray}
     && -2i\tr_{L^2(S)}[\gamma\phi T'(\theta)Ae^{-t\tilde D_\theta^2}]\nonumber\\
     &=&-2i\int_0^\infty \phi(x) \tr_H[\gamma T'(\theta)A e^{-t\tilde D_\theta^2}(x,x)]dx
        \nonumber\\
     &=:&\ti(t)+\tii(t)+\tiii(t).\mylabel{G6-3.34}
\end{eqnarray}
We compute
\begin{equation}
   \ti(t)=-2i(4\pi t)^{-1/2}\int_0^\infty \phi(x)\tr_H[\gamma T'(\theta)A e^{-tA^2}] dx=0,
\end{equation}
since $\gamma$ commutes with $T'(\theta)$ but anticommutes with $A$.
Next we see that
\begin{eqnarray}
   \tii(t)&=&-2 i (4\pi t)^{-1/2}
     \int_0^\infty \phi(x) e^{-x^2/t} \,\tr\left[\gamma T'(\theta) A(I-2 P(\theta))
         e^{-t A^2}\right] dx\nonumber\\
   &\sim_{t\to 0+}& i\,  \tr\left[\gamma T'(\theta) A P(\theta)
         e^{-t A^2}\right].
       \mylabel{G6-3.35}
\end{eqnarray}
Finally, with $d(\gl)=\tr_{\ker(|A|-\gl)} [\gamma T'(\theta) A P(\theta)
         e^{-t A^2}],$

\begin{eqnarray}
   \tiii(t)&\sim_{t\to 0+}&-2 ia(\theta)
    \sum_{\gl\in\spec |A| \setminus\{0\}}
    d(\gl) \sqrt{t}\gl
    \int_0^\infty e^{-2 a(\theta) \gl\sqrt{t}z-t\gl^2}
     \erfc(z)dz\nonumber\\
   &=& -2ia(\theta) \sum_{\gl\in\spec |A|\setminus\{0\}}
       d(\gl) F_{a(\theta)}(\sqrt{t}\gl)\nonumber\\
    &=&-\frac{a(\theta)}{\pi}\int_{\Re w=c}
        t^{-w/2} \eta(A,\gamma T'(\theta) A P(\theta);w-1) {\cm} F_{a(\theta)}(w) dw.
       \mylabel{G6-3.36}
\end{eqnarray}

The existence of the asymptotic expansion hence follows from
our assumptions, Lemma \plref{S6-3.2}, and
\myref{G6-3.32}, \myref{G6-3.35}, \myref{G6-3.36}.
Consequently, we obtain with \myref{G4-1.11a},
\myref{G4-1.10}, \myref{G5-1.23}, and \myref{G7-1.19}:
\begin{eqnarray*}
   a_{-1/2,1}(\check D_\theta, \frac{d}{d\theta} \check D_\theta)&=&
   -(4\pi)^{-1/2} a_{0,1}(A,\gamma i T'(\theta)\\
   &=& \frac{1}{4\sqrt{\pi}} \res(\gamma i T'(\theta)),\\
   a_{-1/2,0}(\check D_\theta, \frac{d}{d\theta} \check D_\theta)&=&
       - (4\pi)^{-1/2} a_{00}(A,\gamma i T'(\theta))\\
       &&+ a_{-1/2,0} (A,\gamma i T'(\theta)A P(\theta))\\
       &&-2a(\theta) {\cm} F_{a(\theta)}(1)
         \Res_1 \eta(A,\gamma i T'(\theta)A P(\theta);0)\\
   &=&- (4\pi)^{-1/2} a_{00}(A,\gamma i T'(\theta))\\
    &&  +\Big(\frac{\sqrt{\pi}}{2}-2a(\theta) {\cm} F_{a(\theta)}(1)\Big)
        \res(\gamma i T'(\theta)(\sgn A) P(\theta)).
\end{eqnarray*}
In view of \myref{G4-1.18} we reach the conclusion.\endproof

\begin{small}
\def\and{{\rm and }}

\end{small}


\begin{thebibliography}{GrSe12}
\bibitem[APS]{APS}
{\sc M.~F.~Atiyah, V.~K.~Patodi \and I.~M.~Singer:}
\newblock {\it Spectral asymmetry and Riemannian geometry I,II,III.}
\newblock Math. Proc. Camb. Phil. Soc., I: {\bf 77} (1975), 43--69;
\newblock II: {\bf 78} (1975), 405--432;
\newblock III: {\bf 79} (1976), 71--99

\bibitem[BL]{BL}
{\sc J.~Br\"uning \and M.~Lesch:}
\newblock {\it On the spectral geometry of algebraic curves.}
\newblock J. reine angew. Math. {\bf 474} (1996), 25--66

\bibitem[BS1]{BS}
{\sc J. Br\"uning \and R. Seeley:}
\newblock {\it The resolvent expansion for second order 
regular singular operators.}
\newblock J. Funct. Anal. {\bf 73} (1987), 369--429

\bibitem[BS2]{BS1}

{\sc J.~Br\"uning \and R.~Seeley:}
\newblock {\it An index theorem for first order regular singular operators.}
\newblock Amer. J. Math. {\bf 110} (1988), 659--714

\bibitem[B]{Bunke}
{\sc U.~Bunke:}
\newblock {\it On the gluing problem for the $\eta$--invariant.}
\newblock J. Diff. Geom. {\bf 41} (1995), 397--448

\bibitem[C]{Calderon}
{\sc A.~P.~Calder{\'o}n:}
\newblock {\it Boundary value problems for elliptic equations.}
\newblock Outlines of the joint Soviet--American symposium on partial
differential equations, Novosibirsk 1963, 303--304



\bibitem[Ch]{Cheeger}
{\sc J.~Cheeger:}
\newblock {\it Analytic torsion and the heat equation.}
\newblock Ann. Math. {\bf 109} (1979), 259-322.



\bibitem[DF]{DF}
{\sc X.~Dai \and D.~S.~Freed:}
\newblock {\it $\eta$--invariants and determinant lines.}
\newblock J. Math. Phys. {\bf 35} (1994), 5155--5194

\bibitem[DG]{DG}
{\sc J.~J.~Duistermaat \and V.~W.~Guillemin:}
\newblock {\it The spectrum of positive elliptic operators and
periodic bicharacteristics.}
\newblock Invent. Math. {\bf 29} (1975), 39--79

\bibitem[DW]{Dougwoj}
{\sc R.~G.~Douglas \and K.~P.~Wojciechowski:}
\newblock {\it 
Adiabatic limits of the $\eta$-invariants: the odd--dimensional
Atiyah--Patodi--Singer problem.}
\newblock Comm. Math. Phys. {\bf 142} (1991), 139--168

\bibitem[GSm]{GilkeySmith}
{\sc P.~B.~Gilkey \and L.~Smith:}
\newblock {\it The eta invariant for a class of elliptic boundary
value problems.}
\newblock Comm. Pure. Appl. Math. {\bf 36} (1983), 85--132

\bibitem[G]{Gi}
{\sc P.~B.~Gilkey:}
\newblock {\it 
Invariance theory, the heat equation, and the Atiyah--Singer index theorem.}
\newblock Second Edition. CRC Press: Boca Raton 1995

\bibitem[GrSe1]{GrubbSeeley0}
{\sc G.~Grubb \and R.~Seeley:},
\newblock{\it Zeta and eta functions for Atiyah--Patodi--Singer operators.}
\newblock To appear in J. Geometric Analysis

\bibitem[GrSe2]{GrubbSeeley}
{\sc G.~Grubb \and R.~Seeley:}
\newblock {\it Weakly parametric pseudodifferential operators and
Atiyah--Patodi--Singer boundary problems.}
\newblock Invent. Math. {\bf 121} (1995), 481--529

\bibitem[HMM]{MazMelHas}
{\sc A.~Hassell, R.~R.~Mazzeo \and R.~B.~Melrose:}
\newblock {\it Analytic surgery and the accumulation of eigenvalues.}
\newblock Comm. Anal. Geom. {\bf 3} (1995), 115--222


\bibitem[Kas]{Kassel}
{\sc C.~Kassel:}
\newblock {\it Le r{\'e}sidue non commutatif (d'apr\`es M. Wodzicki).}
\newblock S\'eminaire Bourbaki, 41\`eme ann\'ee, 1988--89, no 708,
Ast\'erisque {\bf 177--178} (1989), 199--229

\bibitem[Kato]{Kato}
{\sc T.~Kato:}
\newblock {\it Perturbation theory for linear operators.}
\newblock Springer-Verlag: Berlin, Heidelberg, New York 1966

\bibitem[L1]{Lesch1}
{\sc M.~Lesch:}
\newblock {\it 
Deficiency indices for symmetric Dirac operators on manifolds with conical singularities.}
\newblock Topology {\bf 32} (1993), 611--623


\bibitem[L2]{Lesch2}
{\sc M.~Lesch:}
{\it Operators of Fuchs type, conical singularities, 
and asymptotic methods.}
\newblock Teubner--Texte zur Mathematik No. 136, B. G. Teubner Verlag:
Stuttgart--Leipzig, 1997

\bibitem[L3]{Lesch3}
{\sc M.~Lesch:}
\newblock {\it A singular elliptic estimate and applications.}
\newblock In: Pseudo--Differential Calculus and Mathematical
Physics (M. Demuth, E. Schrohe and B.W. Schulze eds.).
Advances in Partial Differential Equations. Akademie Verlag: Berlin 1994,
259--276

\bibitem [LW]{LeschW}
{\sc M.~Lesch \and K.~P.~Wojciechowski:}
\newblock {\it On the $\eta$--invariant of generalized Atiyah--Patodi--Singer
boundary value problems.}
\newblock Ill. J. Math. {\bf 40} (1996), 30--46

\bibitem[MM]{MazMel}
{\sc R.~R.~Mazzeo \and R.~B.~Melrose:}
\newblock {\it Analytic surgery and the eta invariant.}
\newblock Geom. Funct. Anal. {\bf 5} (1995), 14--75

\bibitem[M1]{Mullertorsion}
{\sc W.~M\"uller:}
\newblock {\it Analytic torsion and R--torsion of Riemannian
manifolds.}
\newblock Adv. Math. {\bf 28} (1978), 233--305

\bibitem[M2]{Muller}
{\sc W.~M\"uller:}
\newblock {\it Eta invariants and manifolds with boundary.}
\newblock J. Diff. Geom. {\bf 40} (1994), 311--377

\revmark
\bibitem[M3]{Muller1}
{\sc W.~M\"uller:}
\newblock {\it On the $L^2$--index of Dirac operators on manifolds
with corners of codimension two. I.}
\newblock J. Diff. Geom. {\bf 44} (1996), 97--177

\bibitem[P]{Palais}
{\sc R.~S.~Palais:}
\newblock {\it Seminar on the Atiyah--Singer index theorem.}
\newblock Princeton University Press: Princeton 1965


\bibitem[Si]{Singer}
{\sc I.~M.~Singer:}
\newblock {\it The eta invariant and the index.}
\newblock In: Mathematical aspects of string theory
(S.T. Yau ed.). World Scientific: Singapore 1988, 239--258

\bibitem[So]{Sommer}
{\sc A.~Sommerfeld:}
\newblock {\it Vorlesungen \"uber Theoretische Physik. Bd. 6.}
\newblock Akademische Verlagsgesellschaft: Leipzig 1962

\bibitem[V]{Vishik}
{\sc S.~M.~Vishik:}
\newblock {\it Generalized Ray--Singer Conjecture I. A manifold with
smooth boundary.}
\newblock Commun. Math. Phys. {\bf 167} (1995), 1--102


\bibitem[W]{Witten}
{\sc E.~Witten:}
\newblock {\it Global gravitational anomalies.}
\newblock Commun. Math. Phys. {\bf 100} (1985), 197--229

\revmark
\bibitem[Wod1]{Wod1}
{\sc M.~Wodzicki:}
\newblock {\it  Spectral asymmetry and local invariants (in Russian).}
\newblock Habilitation thesis, Moscow: Steklov Math. Inst., 1984

\revmark
\bibitem[Wod2]{Wod2}
{\sc M.~Wodzicki:}
\newblock {\it Non--commutative residue I.}
\newblock Lecture Notes in Mathematics {\bf 1289} (1987), 320--399

\bibitem[W1]{Woj1}
{\sc K.~P. Wojciechowski:}
\newblock {\it The additivity of the $\eta$--invariant. The case of
an invertible tangential operator.}
\newblock Houston J. Math. {\bf 20} (1994), 603--621

\bibitem[W2]{Woj2}
{\sc K.~P. Wojciechowski:}
\newblock {\it The additivity of the $\eta$--invariant. The case of
a singular tangential operator.}
\newblock Commun. Math. Phys. {\bf 109} (1995), 315--327 


\end{thebibliography}
\end{document}